\documentclass[moor,nonblindrev]{informs1} % for review, not blinded
\usepackage{natbib}
 \NatBibNumeric
 \def\bibfont{\small}%
 \bibpunct[, ]{[}{]}{,}{n}{}{,}%

\usepackage[colorlinks=true,breaklinks=true,bookmarks=true,urlcolor=blue,
     citecolor=blue,linkcolor=blue,bookmarksopen=false,draft=false]{hyperref}

\def\EMAIL#1{\href{mailto:#1}{#1}}% When hyperref is used, otherwise outcomment 
\def\URL#1{\href{#1}{#1}}         % When hyperref is used, otherwise outcomment 

\TheoremsNumberedThrough     % Preferred (Theorem 1, Lemma 1, Theorem 2)
\EquationsNumberedThrough    % Default: (1), (2), ...
\usepackage[normalem]{ulem}

\usepackage[]{algorithm}
\usepackage{algorithmic}              

\algsetup{indent=12pt}

\usepackage{tikz}
\usetikzlibrary{fit,positioning,arrows,automata}   

\usepackage{float}
\usepackage{subfigure}        
                 
\usepackage{mathtools}           
\mathtoolsset{showonlyrefs}

\def\minwrt[#1]{\underset{#1}{\minimize}\;}
\def\argminwrt[#1]{\underset{#1}{\text{arg min }}}
\def\maxwrt[#1]{\underset{#1}{\maximize}\;}
\def\argmaxwrt[#1]{\underset{#1}{\text{arg max }}}
\def\maxemphwrt[#1]{\underset{#1}{\text{\emph{maximize} }}}

\def\bK{{\bf K}}
\def\bM{{\bf M}}

\def\bU{{\bf U}}

\def\bC{{\bf C}}

\def\tr{{\rm trace}}
\def\ett{{\bf 1}}
\def\ccE{{\mathcal{E}}}

\def\ccH{{\mathcal{H}}}
\def\ccL{{\mathcal{L}}}

\def\ccS{{\mathcal{S}}}
\def\ccT{{\mathcal{T}}}
\def\ccO{{\mathcal{O}}}
\def\ccP{{\mathcal{P}}}

\def\ccV{{\mathcal{V}}}

\def\RR{{\mathbb{R}}}

\def\RRext{{\overline{\mathbb{R}}}}

\def\Nflow{{\mathcal{N}}}
\def\Vflow{{\mathcal{V}}}
\def\Eflow{{\mathcal{E}}}

\def\Gomt{{G}}
\def\Vomt{{V}}
\def\Eomt{{E}}

\def\CL{ {C_L}}
\def\KL{ {K_L}}

\DeclareMathOperator{\minimize}{minimize}
\DeclareMathOperator{\maximize}{maximize}

\newcommand{\diag}{{{\mathrm{diag}}}}

\begin{document}
%%%%%%%%%%%%%%%%

%\setcounter{page}{1}
%\DOI{10.1287/xxxx.0000.0000}%

\RUNAUTHOR{Haasler et al.}

\RUNTITLE{Dynamic flow problems via optimal transport}

% \TITLE{Bundling Information Goods of Decreasing Value}
\TITLE{Scalable computation of dynamic flow problems via multi-marginal graph-structured optimal transport}

\ARTICLEAUTHORS{%
\AUTHOR{Isabel Haasler}
\AFF{KTH Royal Institute of Technology, \EMAIL{haasler@kth.se}, \URL{}}
\AUTHOR{Axel Ringh}
\AFF{The Hong Kong University of Science and Technology, \EMAIL{eeringh@ust.hk}, \URL{}}
\AUTHOR{Yongxin Chen}
\AFF{Georgia Institute of Technology, \EMAIL{yongchen@gatech.edu}, \URL{}}
\AUTHOR{Johan Karlsson}
\AFF{KTH Royal Institute of Technology, \EMAIL{johan.karlsson@math.kth.se}, \URL{}}
} % end of the block

\ABSTRACT{%
 % Enter your abstract

In this work, we develop a new framework for dynamic network flow problems based on optimal transport theory.
We show that the dynamic multi-commodity minimum-cost network flow problem can be formulated as a multi-marginal optimal transport problem, where the cost function and the constraints on the marginals are associated with a graph structure. 
By exploiting these structures and building on recent advances in optimal transport theory, we develop an efficient method for such entropy-regularized optimal transport problems.
In particular, the graph structure is utilized to efficiently compute the projections needed in the corresponding Sinkhorn iterations, and we arrive at a scheme that is both highly computationally efficient and easy to implement.
To illustrate the performance of our algorithm, we compare it with a state-of-the-art Linear programming (LP) solver. We achieve good approximations to the solution at least one order of magnitude faster than the LP solver.
Finally, we showcase the methodology on a traffic routing problem with a large number of commodities.

}%

\KEYWORDS{Multi-marginal optimal transport; Dynamic network flow; Multi-commodity network flow; Sinkhorn’s method; Computational methods; Traffic routing}
\MSCCLASS{49Q22, 90-08, 90C35, 90C08, 49M29 }
\ORMSCLASS{
Networks/graphs: Flow algorithms; Networks/graphs: Multicommodity; Programming: Linear: Large scale systems} %Primary: ; secondary: }
\HISTORY{}

\maketitle
%%%%%%%%%%%%%%%%%%%%%%%%%%%%%%%%%%%%%%%%%%%%%%%%%%%%%%%%%%%%%%%%%%%%%%

% Samples of sectioning (and labeling) in MOOR.
% NOTE: (1) all section levels end with a period,
%       (2) capitalization is as shown (sentence style, not title style).
%
%\section{Introduction.}\label{intro} %%1.
%\subsection{Duality and the classical EOQ problem.}\label{class-EOQ} %% 1.1.
%\subsection{Outline.}\label{outline1} %% 1.2.
%\subsubsection{Cyclic schedules for the general deterministic SMDP.}
%  \label{cyclic-schedules} %% 1.2.1
%\section{Problem description.}\label{problemdescription} %% 2.

% Text of your paper here

\section{Introduction.}

Many phenomena in today's society can be modelled as large scale transportation or flow problems, and new technological advances create the need for solving larger and larger problems. An example is the introduction of self driving-cars to the road network, which will create both new opportunities and new challenges \cite{levinson2011towards, pasquale2019traffic}.
Increasing automation and communication between vehicles will result in very large systems where all vehicles need to be routed simultaneously taking into account destinations, vehicle properties and urgency \cite{carlino2012approximately}.
Another challenge is to direct large crowds in, e.g., transit areas in airports, subways, or event venues \cite{yamada1996network, haghani1996formulation, aronson1989survey}, which is particularly critical for evacuation scenarios in the case of emergencies, but also essential for every-day use.

Many of these problems can be modelled as large scale dynamic network flow problems \cite{bertsimas2000traffic, kennington1978survey, aronson1989survey}.
The most common strategy for handling such problems is to convert the dynamic flow problem to a static 
flow problem on a time-expanded network, and this strategy goes back to the classical work \cite{ford1958constructing}. In addition to this, there are typically several classes of groups of agents with heterogeneous properties and objectives in the system.
For instance, each agent in a traffic network drives a vehicle with certain properties, and the objective is typically to reach a certain destination with a certain degree of urgency. 
Similar problems appear in air traffic planning, railroad traffic scheduling, communication and logistics, and are often treated as multi-commodity flow problems over networks \cite{haghani1996formulation, bertsimas2000traffic, kennington1978survey, aronson1989survey}.
Although such problems are usually formulated as linear programming (LP) problems, for real applications the corresponding optimization problems are often too large to be handled by standard methods.  
Specialized methods exploit the structure of multi-commodity flow problems, using, e.g., column generation methods. These include price-directive decomposition \cite{jones1993multicommodity}, resource-directive decomposition \cite{kennington1977effective, mcbride1998progress}, and basis partitioning methods \cite{farvolden1993primal}.
However, it has been reported that these methods typically decrease the solution time of standard (LP) solvers by at most one order of magnitude \cite{barnhart2009multicommodity, retvdri2004novel}.

During the last few decades there has been considerable development in the field of optimal transport theory.
Traditionally the optimal transport problem addresses a static scenario where one given distribution is transported to another, and this problem has been extensively used in areas such as economics and logistics \cite{villani2008optimal}.
There has recently been a rapid advancement of theory and applications for optimal transport, in particular towards applications in imaging, statistics and machine learning (see \cite{peyre2019computational} and references therein), and systems and control  \cite{benamou2000computational,chen2016optimalPartI}, which has led to a mature framework with computationally efficient algorithms
\cite{peyre2019computational}
that can be used to address a wide range of problems. The optimal transport problem is a linear program, but the number of variables often makes it intractable to solve with general-purpose optimization methods for large size problems. 
However, a recent computational breakthrough in this area builds on introducing an entropic barrier term in the objective function. The resulting optimization problem can then be solved efficiently using the so called Sinkhorn iterations \cite{cuturi2013sinkhorn}.
This allows for computing an approximate solution of large transportation problems and has opened up the field for new applications where no computationally feasible method previously existed.

The optimal transport framework has in some cases been used for modelling several kinds of interacting classes, e.g., for transport of multiple species \cite{chen2018vector,bacon2020multi} or flows with several phases \cite{benamou2004numerical}.
In this paper we will build on some of these results and we propose to use a generalization of the optimal transport problem with several marginals to address multi-commodity flow problems.
This multi-marginal optimal transport problem
\cite{gangbo1998optimal, pass2015multi, ruschendorf1995optimal, ruschendorf2002n} is computationally challenging since the number of variables grows exponentially in the number of marginals.
Even though entropy regularization methods have been derived for the multi-marginal optimal transport problem \cite{benamou2015bregman},
the cost for each iteration still grows exponentially in the number of marginals (see \cite{lin2019complexity} for computational complexity bounds).
However, in many cases the cost function has a structure that can be utilized for efficient computations, as for example in barycenter, information fusion, and tracking problems \cite{benamou2015bregman, elvander19multi,haasler20trees}.

In this paper we show that the dynamic flow problem can be formulated as a structured multi-marginal optimal transport problem.
This structure can be visualized in a graph where the set of nodes corresponds to the marginals, and where there is an edge between two nodes if there is a cost term or a constraint that depends jointly on the two nodes.
For the single commodity case, this structure is a path graph with one node for each time point that represents the flow in the network at that time.
For the dynamic multi-commodity network flow problem, there is one additional node in the graph that represents the distribution over the different commodity classes.
The solution to this optimal transport problem then describes a joint distribution, which consists of the optimal flow for all commodities in the dynamic network problem.

We consider the corresponding entropy-regularized approximation of this problem, and by utilizing the structure in the cost function we derive methods for solving this problem.
Many of the classical methods for dynamic flow problems consider standard network flow methods on the time-expanded network. By instead formulating this problem as a multi-marginal optimal transport problem, we can more efficiently utilize the sequential structure without explicitly setting up the time-expanded network. This results in an elegant and easily implementable method.
We illustrate experimentally that this method is computationally competitive with state-of-the-art methods, and then apply it to a traffic routing problem.

The rest of the paper is structured as follows.
Section~\ref{sec:background} summarizes background material on dynamic multi-commodity network flows and multi-marginal optimal transport.
In Section~\ref{sec:nf_ot} we explain how to formulate network flow problems as structured multi-marginal optimal transport problems. Based on this, we develop numerical schemes to solve the problems in Section~\ref{sec:OTgraph}.
Finally, in Section~\ref{sec:siumlations} we compare the performance of our methods to a commercial LP solver, and showcase it in a traffic routing application.

\section{Background.} \label{sec:background}

In this section we review background on the two central topics of this paper: dynamic multi-commodity network flows and multi-marginal optimal transport.
We also use this Section to set up notation. In particular, bold-faced letters are used throughout to denote tensors, and $\otimes$ denotes the tensor (outer) product, e.g., for vectors $v_1 \in \RR^{n_1}$ and $v_2 \in \RR^{n_2}$ we have that $v_1 \otimes v_2 \in \RR^{n_1 \times n_2}$ and $(v_1 \otimes v_2)_{ij}  = (v_1)_i(v_2)_j$. Moreover, by $\ett$ we denote a column vector of ones of appropriate size, by $\mathbb{R}_+$ we denote the nonnegative real numbers, and we use $\RRext_+= \RR_+ \cup \{ \infty \}$ and $\RRext= \RR \cup \{ \infty \}\cup\{ -\infty \}$ to denote the extended nonnegative real line and extended  real line, respectively. Throughout we will adopt the convention that $0\cdot \infty=0$. 
Finally, by $\exp(\cdot)$, $\log(\cdot)$, $\odot$, $./$, and $\min(\cdot, \cdot)$ we denote elementwise exponential, logarithm, product, division, and minimum respectively.

%\section{Minimum-cost network flow problems}
\subsection{Minimum-cost network flow problems.}
\label{sec:network_flow} 

A minimum-cost network flow problem is to determine a flow from sources to sinks with minimum cost  \cite{ford1962flows, bertsekas1988relaxation}. More specifically, the flow is defined on a network $\Nflow=(\Vflow,\Eflow)$ with vertices $\Vflow$ and directed edges $\Eflow$, and the sources and sinks are sets of edges%
\footnote{Often the sources and sinks are defined on the nodes $\Vflow$ not the edges $\Eflow$. In this work we consider the latter case, however the framework introduced herein can easily be modified to define the sources and sinks on the nodes $\Vflow$ instead.} 
$\ccS^+ \subset \Eflow$ and $\ccS^- \subset \Eflow$.
Let each source $e\in \ccS^+$ be equipped with a supply $r^+_e\in\RR_{+}$, and each sink $e\in \ccS^-$ with a demand $r^-_e\in\RR_{+}$, and we assume that the total supply matches the total demand, i.e., that $\sum_{e\in \ccS^+} r^+_e - \sum_{e\in \ccS^-} r^-_e = 0$.
In addition, let each edge $e\in \Eflow$ be assigned a cost $c_e\in \RR_+$ of transporting a unit of flow on that edge. The goal of minimum cost-flow problems is to transport the flow from the sources to the sinks with minimal total transporting cost. 
We also include capacity constraints, which require that the total flow on an edge is limited by the edge capacity $d_e \in\RR_{+}$ on $e\in\Eflow$.

There are two standard formulations for the network flow problem. 
One is the arc-chain formulation, where one optimizes over a set of flow paths (arc-chains) from sources to sinks \cite{ford1962flows, tomlin1966minimum}. This is the main formulation considered in this work and is described in detail below. Another common formulation is the node-edge formulation, where one seeks the optimal amount of flow over each edge while maintaining flow balance in each node. 
For more details on this formulation, and a comparison of both formulations we refer the reader to \cite{ford1962flows, tomlin1966minimum}.

\subsubsection{The arc-chain formulation.}
%\subsection{The arc-chain formulation}
Given a network $\Nflow=(\Vflow,\Eflow)$, a path is a sequence of edges that joins two vertices such that all edges and all visited vertices are distinct, i.e., they occur at most once in the sequence \cite[p.~6]{diestel2017graph}.
A path is thus a subgraph, which we denote by $p$, and is defined by a list of edges $(p_1,p_2,\dots,p_N)$, where $ p_j \in \ccE$  denotes the $j$-th element of the path for $j=1,\dots,N$. Here, $N$ is called the length of the path $p$. Moreover, since $p$ is a path the edge $p_j$ ends in the initial node of $p_{j+1}$ for $j=1,\dots,N-1$.

In the arc-chain formulation, we consider the paths, or arc-chains, which start in a source and end in a sink.
Let $\ccP$ denote the set of all such paths, where the first element lies in $\ccS^+$, and its last element lies in $\ccS^-$.
Moreover, let $\ccP^+_{e}$ denote the paths starting from the edge $e \in \ccS^+$, and let $\ccP^-_{e}$ denote the paths ending in the edge $e \in \ccS^-$.
The cost of a path $p\in \ccP$ is the sum of the costs of its edges  $c_p = \sum_{e\in p}  c_e$.
Next, let $x_p$ denote the amount of flow associated with path $p\in \ccP$.
Then, the arc-chain formulation of the minimum-cost network flow problem reads
\begin{equation} \label{eq:flow_single-commodity}
	\begin{aligned}
\minwrt[x_p \in \RR_+,\,p \in \ccP] &  \sum_{p \in \ccP}  c_{p} x_p\\
\text{subject to} & \sum_{ p \in \ccP^+_e } x_p = r^+_e, \quad \text{ for } e \in \ccS^+, \\
& \sum_{ p \in \ccP^-_e } x_p = r^-_e, \quad \text{ for } e \in \ccS^-,\\ 
&  \sum_{p \in \ccP} \delta_{e \in p} x_p \leq d_e, \quad \text{ for } e\in \Eflow,
\end{aligned}
\end{equation}
where $\delta_{e \in p} = 1$ if the edge $e$ is part of path p, and $\delta_{e \in p} = 0$ otherwise.
Here, the objective function corresponds to the total cost of the flow.
The first two sets of constraints guarantee that the supply and demand for all sources and sinks are satisfied, and the last set of constraints enforces that the flow on each edge does not exceed the given capacity.

\subsubsection{Multi-commodity network flow.}
%\subsection{Multi-commodity network flow}

The extension to multi-commodity network flow problems deals with the case where there are multiple commodities  present in the network \cite{hall2007multicommodity, wang2018multicommodity, tomlin1966minimum, ford1958suggested, kennington1978survey}. Here we let $L$ denote the number of commodities, and let $c_{e}^\ell$ denote the cost of a unit flow on edge $e\in \Eflow$ of commodity $\ell$, for $\ell=1,\dots,L$. The supply and demand typically depend on the commodity, thus 
each commodity $\ell$ has specified sources $\ccS^{\ell,+} \in \ccE$ with supplies $r_e^{\ell,+}$ for $e\in \ccS^{\ell,+} $, and sinks $\ccS^{\ell,-} \in \Eflow$  with demands  $r_e^{\ell,-}$ for $e\in \ccS^{\ell,-} $.
Moreover, for each commodity $\ell=1,\dots,L$, let $\ccP^\ell$ denote the sets of paths from the sources to the sinks, and let $\ccP^{\ell,+}_{e}$ denote the paths starting in $e \in \ccS^{\ell,+}$, and let $\ccP^{\ell,-}_{e}$ denote the paths ending in $e \in \ccS^{\ell,-}$.
The cost of a unit flow of commodity $\ell$ on a path $p\in \ccP$ is the sum of the corresponding costs of the edges in the path $c_p^\ell = \sum_{e\in p}  c_e^\ell$.
Next, by letting $x_p^\ell$ denote 
 the amount of flow of commodity $\ell$ on path $p$, the minimum cost multi-commodity network flow problem in arc-chain formulation reads
\begin{equation} \label{eq:flow_multi-commodity}
\begin{aligned}
\minwrt[\substack{x_p^\ell \in \RR_+, \, p \in \ccP^\ell \\\ell =1,\dots, L}] & \sum_{\ell=1}^{L} \sum_{p \in \ccP^\ell}  c_{p}^\ell x_p^\ell\\
\text{subject to} & \sum_{p \in \ccP^{\ell,+}_e} x_p^\ell = r^{\ell,+}_e ,  \quad \mbox{  for }  e \in \ccS^{\ell,+}, \quad  \ell =1,\dots, L, 
\\  
&\sum_{p \in \ccP^{\ell,-}_e} x_p^\ell = r^{\ell,-}_e , \quad   \mbox{  for } e \in \ccS^{\ell,-}, \quad  \ell =1,\dots, L, \\
& \sum_{\ell=1}^{L} \sum_{p \in \ccP^\ell} \delta_{e \in p} x_p^\ell \leq d_e, \quad \text{ for } e\in \Eflow.
\end{aligned}
\end{equation}
Here, the first two sets of constraints guarantee that the demand and supply for all commodities are satisfied. The third set of constraints enforces that the flow on each edge does not exceed the given capacity. In particular, note that the multi-commodity problem \eqref{eq:flow_multi-commodity} with only one commodity, i.e., $L=1$, boils down to the single-commodity problem \eqref{eq:flow_single-commodity}.

\subsubsection{Dynamic network flow.}
%\subsection{Dynamic network flow}

In this work we consider dynamic flows, also called flows over time, where the time that it takes for the flow to travel in the network is taken into account  \cite{ford1958constructing, aronson1989survey, hall2007multicommodity}. In this work we develop efficient methods that exploit the temporal structure. For this to work we need to assume synchronous travelling times for all edges, but on the other hand the efficient methods  allows for handling problems with large networks and fine time discretization.

More precisely, we consider a flow problem on the network $\Nflow=(\ccV, \ccE)$ over the time interval $0$ to $\ccT$.  The problem is to transport a given flow at time $0$  through the network to a final flow at time $\ccT$ with minimal cost, while satisfying capacity constraints at all time points.  We consider the discretized problem on the time steps $0,1,\ldots, \ccT$.  Dynamic flow problems are typically solved as a static problem on the time-expanded network \cite{ford1958constructing}.
The time-expanded network $\Nflow_{\rm exp}$ is constructed by considering $\ccT+1$ copies of the vertices $\Vflow$, denoted by $\Vflow_0,\dots,\Vflow_\ccT$. Here the copy $\Vflow_t$ is associated with time instance $t$ in the time expanded network, and we denote these nodes by $(t,v)$ where $v\in \Vflow$ in the original network.

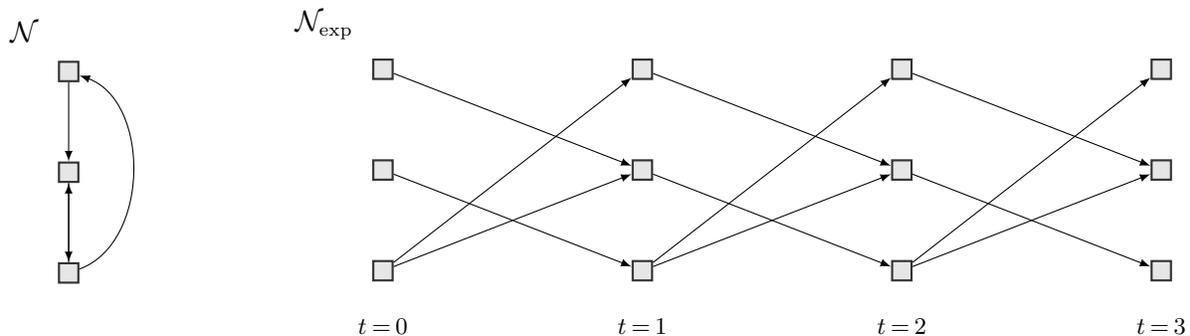
\begin{figure}
\centering
\begin{subfigure}
	\centering
	\begin{tikzpicture}
	\tikzstyle{main}=[rectangle, minimum size = 5pt, thick, draw =black!80, node distance = 30pt and 100pt]
	
	\node[main,fill=black!10] (t0n1) {}; 	
	\node[main,fill=black!10] (t0n2) [below =of t0n1] {};  
	\node[main,fill=black!10] (t0n3) [below =of t0n2] {};
	\node[draw=none] (t0) [below =15pt of t0n3] {};
	
	\node[draw=none] (N) [above left =5pt of t0n1] {$\Nflow$};

	\draw[ -latex] (t0n1) --  (t0n2);
	\draw[ -latex] (t0n2) -- (t0n3);
	\draw[ -latex] (t0n3) to [out=20,in=340] (t0n1); 
	\draw[ -latex] (t0n3) -- (t0n2);
	
	\end{tikzpicture} \hspace{40pt}
\end{subfigure}
\begin{subfigure}
	\centering
	\begin{tikzpicture}
	\tikzstyle{main}=[rectangle, minimum size = 5pt, thick, draw =black!80, node distance = 30pt and 90pt]
	
	\node[main,fill=black!10] (t0n1) {}; 	
	\node[main,fill=black!10] (t0n2) [below =of t0n1] {};  
	\node[main,fill=black!10] (t0n3) [below =of t0n2] {}; 
	\node[draw=none] (t0) [below =10pt of t0n3] {	\footnotesize $t=0$}; 
	\node[draw=none] (Nexp) [above left =5pt of t0n1] {$\Nflow_{\rm exp}$}; 	
 
	\node[main,fill=black!10] (t1n1) [right =of t0n1] {}; 	
	\node[main,fill=black!10] (t1n2) [below =of t1n1] {};  
	\node[main,fill=black!10] (t1n3) [below =of t1n2] {};
	\node[draw=none] (t1) [below =10pt of t1n3] {\footnotesize $t=1$}; 
		
 	\node[main,fill=black!10] (t2n1) [right =of t1n1] {}; 	
	\node[main,fill=black!10] (t2n2) [below =of t2n1] {};  
	\node[main,fill=black!10] (t2n3) [below =of t2n2] {}; 
	\node[draw=none] (t2) [below =10pt of t2n3] {\footnotesize $t=2$}; 	
	
	\node[main,fill=black!10] (t3n1) [right =of t2n1] {}; 	
	\node[main,fill=black!10] (t3n2) [below =of t3n1] {};  
	\node[main,fill=black!10] (t3n3) [below =of t3n2] {}; 
	\node[draw=none] (t3) [below =10pt of t3n3] {\footnotesize $t=3$}; 	
	
	\draw[ -latex] (t0n1) --  (t1n2);
	\draw[ -latex] (t0n2) -- (t1n3);
	\draw[ -latex] (t0n3) -- (t1n1); 
	\draw[ -latex] (t0n3) -- (t1n2);
	
	\draw[ -latex] (t1n1) --  (t2n2);
	\draw[ -latex] (t1n2) -- (t2n3);
	\draw[ -latex] (t1n3) -- (t2n1); 
	\draw[ -latex] (t1n3) -- (t2n2); 
	
	\draw[ -latex] (t2n1) --  (t3n2);
	\draw[ -latex] (t2n2) -- (t3n3);
	\draw[ -latex] (t2n3) -- (t3n1); 
	\draw[ -latex] (t2n3) -- (t3n2); 

%	\draw (52.5pt,-37.5pt) ellipse (30pt and 70pt);
%	\draw (160pt,-37.5pt) ellipse (30pt and 70pt);
%	\draw (267.5pt,-37.5pt) ellipse (30pt and 70pt);

	\end{tikzpicture}
	\end{subfigure}
	\caption{A network with three nodes and its time-expanded network for $\ccT=3$ time steps.} \label{fig:time_exp}
\end{figure}

The edges of $\Nflow_{\rm exp}$ connect nodes corresponding to consecutive time instances according to the edges $\Eflow$ in the original network, that is, 
$\Eflow_{\rm exp}=\cup_{t=1}^\ccT  \Eflow_t$
where 
 $\Eflow_t$ consists of the directed edges $((t-1, v_{t-1}), (t, v_t))$ where 
$(v_{t-1},  v_t)\in \ccE$, for $t=1,\ldots, \ccT$.
The capacities and costs on these added edges are defined to be the same as the corresponding\footnote{Note that there is a canonical bijection $(v_1, v_2) \leftrightarrow ((t-1, v_1), (t, v_2))$ between the edges $\ccE$ and the edges $\ccE_t$.} edges in the original network $\Nflow$.
The time-expanded network is illustrated for a simple example in Figure~\ref{fig:time_exp}.

To express the dynamic flow problem in arc-chain formulations similarly to \eqref{eq:flow_single-commodity} and \eqref{eq:flow_multi-commodity},
a path $p$ is as before a tuple of edges
 $(p_1,\dots,p_\ccT)$. Its element $p_t \in \Eflow_t$ denotes the edge, which the paths flow takes in the time interval  $[t-1,t]$.
In the setting of one commodity, let $\ccP$ denote the set of feasible paths in the time-expanded network $\Nflow_{\rm exp}$, i.e., $p\in\ccP$ if $p$ is a path that starts in a source, $p_1\in \ccS^{+}$, and ends in a sink, $p_\ccT\in \ccS^{-}$.
The corresponding cost of unit flow on the path $p\in \ccP$ is then $c_p=\sum_{t=1}^\ccT c_{p_t}$.
The dynamic minimum-cost network flow problem can then be written as
\begin{subequations} \label{eq:flow_dynamic_single-commodity}
\begin{align}
\minwrt[x_p \in \RR_+ , \, p \in \ccP] &   \sum_{p \in \ccP}   c_p x_p \label{eq:flow_dynamic_single-commodity_objective}\\
\text{subject to} & \sum_{p \in \ccP}  \delta_{e=p_1} x_p = r^{+}_e, \quad \text{ for } e \in \ccS^{+},  \label{eq:flow_dynamic_single-commodity_s} \\
&  \sum_{p \in \ccP} \delta_{e=p_\ccT}x_p = r^{-}_e, \quad \text{ for }  e \in \ccS^{-},   \label{eq:flow_dynamic_single-commodity_d} \\
&  \sum_{p \in \ccP} \delta_{e = p_t} x_p \leq d_e, \quad \text{ for }  e\in \Eflow, \quad t=2,\dots, \ccT-1 . \label{eq:flow_dynamic_single-commodity_capacity}
\end{align}
\end{subequations}
Note that the network flow problem \eqref{eq:flow_single-commodity} on the time-expanded network $\Nflow_{\rm exp}$,  corresponds to \eqref{eq:flow_dynamic_single-commodity} line by line.
 
To formulate the multi-commodity counterpart of the dynamic flow problem \eqref{eq:flow_dynamic_single-commodity}, let $\ccP^\ell$ denote the set of feasible paths in the time-expanded network $\Nflow_{\rm exp}$ for commodity $\ell=1,\dots,L$. 
The corresponding cost of unit flow on the path for commodity $\ell$ is then
$c_p^\ell=\sum_{t=1}^\ccT c_{p_t}^{\ell}$ for a path $p\in \ccP^\ell$, and the dynamic minimum-cost multi-commodity network flow problem reads
\begin{equation} \label{eq:flow_dynamic_multi-commodity}
\begin{aligned}
\minwrt[\substack{x_p^\ell \in \RR_+, \, p \in \ccP^\ell\\ \ell=1,\dots, L}] &  \sum_{\ell=1}^{L}  \sum_{p \in \ccP^\ell}   c_p^\ell x_p^\ell  \\
\text{subject to} & \sum_{p \in \ccP^\ell}  \delta_{e=p_1} x_p^\ell = r^{\ell,+}_e, \quad \text{ for }  e \in \ccS^{\ell,+},  \quad \ell =1,\dots, L, \\
&  \sum_{p \in \ccP^\ell} \delta_{e=p_\ccT}x_p^\ell = r^{\ell,-}_e, \quad \text{ for }  e \in \ccS^{\ell,-}, \qquad  \ell =1,\dots, L,    \\
& \sum_{\ell=1}^{L}  \sum_{p \in \ccP^\ell} \delta_{e = p_t} x_p^\ell \leq d_e, \quad \text{ for }  e\in \Eflow, \quad t=2,\dots, \ccT-1 ;
\end{aligned}
\end{equation}
see \cite{khodayifar2019minimum} for a similar problem formulation.

A problem with the arc-chain formulations is that the number of variables, corresponding to possible paths, grows exponentially with $\ccT$. Thus, standard linear programming methods are not applicable when $\ccT$ is large. A way to circumvent this issue is to use specialized solvers building on, e.g., column generation, or to instead consider the corresponding node-edge formulations of the problem (cf. \cite{ford1962flows, tomlin1966minimum}). In this work we take a different approach that builds on formulating the problem as an optimal transport problem that utilize the structure in the arc-chain formulation.

\subsection{Optimal transport.}
%\section{Optimal transport}
 \label{sec:omt}

 The optimal transport problem is to find a mapping that moves the mass from one distribution to another with minimal cost, based on an underlying metric \cite{villani2008optimal}.
 In this paper we consider the discrete setting, where the two distributions are represented by two non-negative vectors $\mu_1 \in \RR_{+}^{n_1}$, $\mu_2 \in \RR_{+}^{n_2}$ with equal mass. 
In this setting the transport cost is defined in terms of a underlying non-negative cost matrix
$C \in \RRext_{+}^{n_1\times n_2}$, where $C_{ij}$ denotes the cost%
\footnote{If transport of mass is not allowed from position $i$ to position $j$, then we let $C_{ij}=\infty$.}
of moving a unit mass from position $i$ to $j$.
Analogously, a transport plan $M \in \RR_+^{n_1\times n_2}$ is a non-negative matrix, where $M_{ij}$ represents the amount of mass moved from $i$ to $j$.
The optimal transport plan from $\mu_1$ to $\mu_2$ is then a minimizing solution of
\begin{equation} \label{eq:omt}
\begin{aligned}
\minwrt[M \in \RR_+^{n_1\times n_2}] \ &\ \tr (C^T M) \\
\text{ subject to } \ &\ M \ett = \mu_1,\\
&\  M^T \ett = \mu_2 .
\end{aligned}
\end{equation}
%
%\subsection{Multi-marginal optimal transport}
%\subsubsection{Multi-marginal optimal transport}
Multi-marginal optimal transport extends the concept of the classical optimal transport problem \eqref{eq:omt} to the setting with a set of marginals $\mu_t \in \RR_{+}^{n_t}$, for $t=1,\dots,\ccT$,
%$\mu_1,\dots,\mu_\ccT \in \RR^n$, 
where $\ccT\geq 2$ \cite{pass2015multi,benamou2015bregman,elvander19multi, haasler20trees}. In this setting, the transport cost and transport plan are described by tensors $\bC \in \RRext^{n_1\times n_2 \dots \times n_\ccT}_+$ and $\bM \in \RR^{n_1\times n_2 \dots \times n_\ccT}_+$. Here, $\bC_{i_1\dots i_\ccT}$ denotes the unit cost associated with the tuple $(i_1,\dots,i_\ccT)$, and $\bM_{i_1\dots i_\ccT}$ denotes the amount of mass associated with this tuple.
Then the total transportation cost for a given transport plan $\bM$ is
\begin{equation}
\langle \bC, \bM \rangle = \sum_{i_1,\dots,i_\ccT} \bC_{i_1 \dots i_\ccT} \bM_{i_1 \dots i_\ccT}.
\end{equation}
Moreover, $\bM$ is a transport plan between the desired marginals if its projections on the marginals satisfy $P_t(\bM) = \mu_t$, for $t=1,\dots,\ccT$, where the projection on the $t$-th marginal is defined by
\begin{equation} \label{eq:proj_discrete}
(P_t(\bM))_{i_t} := \sum_{i_1,\dots,i_{t-1},i_{t+1},\ldots, i_\ccT} \bM_{i_1\dots i_{t-1} i_t i_{t+1} \dots i_\ccT}.
\end{equation}
The discrete multi-marginal optimal transport problem thus reads
\begin{equation} \label{eq:omt_multi_discrete}
\begin{aligned}
\minwrt[ \bM \in \RR^{n_1 \times \dots \times n_\ccT}_+] \ &\ \langle \bC, \bM \rangle  \\
\text{ subject to } \ &\ P_t (\bM) = \mu_t,  \text { for } t \in \Gamma.
\end{aligned}
\end{equation}
 Here $\Gamma$ is an index set that describes the set of constrained marginals. In the original multi-marginal optimal transport formulation, constraints are typically given on all marginals, i.e., for the index set $\Gamma = \{1,2,\dots, \ccT \}$. However, in this work we typically consider the case where constraints are only imposed on a subset of marginals, i.e., $\Gamma \subset \{1,2,\dots, \ccT\}$, or when some of the constraints are inequality constraints.

Note that the standard bi-marginal optimal transport problem \eqref{eq:omt} is a special case of the multi-marginal optimal transport problem \eqref{eq:omt_multi_discrete}, where $\ccT=2$ and $\Gamma=\{1,2\}$.
It is also worth noting that the bi-marginal optimal transport problem can be interpreted as a minimum-cost network flow problem. However, this interpretation does in general not extend to the multi-marginal case \cite{lin2020fixed}.
In this work we show how to formulate any dynamic network flow problem as a multi-marginal optimal transport problem with a structured cost tensor.

%\subsection{Sinkhorn iterations}
\subsubsection{Sinkhorn iterations.}
 \label{subsec:sinkhorn}
Although linear, the number of variables in the multi-marginal optimal transport problem \eqref{eq:omt_multi_discrete} is often too large to be solved directly. 
A popular approach for the bi-marginal setting to bypass the size of the problem has been to add a regularizing entropy term to the objective \cite{cuturi2013sinkhorn}. In principle, the same approach can be used also for the multi-marginal case. 
With the entropy term
\begin{equation} \label{eq:entropy_term} 
D(\bM) = \sum_{i_1,\dots,i_\ccT} \left( \bM_{i_1\dots i_\ccT} \log(\bM_{i_1\dots i_\ccT}) + \bM_{i_1\dots i_\ccT} -1 \right),
\end{equation}
the entropy regularized multi-marginal optimal transport problem is defined as
\begin{equation} \label{eq:omt_multi_regularized}
\begin{aligned}
\minwrt[ \bM \in \RR^{n_1\times \dots \times n_\ccT}_+] & \langle \bC, \bM \rangle + \epsilon D(\bM) \\
\text{ subject to } & P_t (\bM) = \mu_t,  \text { for } t \in \Gamma,
\end{aligned}
\end{equation}
where $\epsilon>0$ is a small regularization parameter.
The introduction of the entropy term in problem \eqref{eq:omt_multi_regularized} allows for expressing the optimal solution $\bM$ in terms of Lagrange dual variables, which may be computed by Sinkhorn iterations \cite{benamou2015bregman, nenna2016numerical}.
In particular, it can be shown that the optimal solution of \eqref{eq:omt_multi_regularized} is of the form \cite{elvander19multi}
\begin{equation} \label{eq:MKU}
\bM = \bK \odot \bU ,
\end{equation}
where $\bK = \exp(- \bC/\epsilon)$ and where $\bU$ can be decomposed as
\begin{equation}\label{eq:U}
\bU= u_1 \otimes u_2 \otimes \dots \otimes u_\ccT .
\end{equation}
Here, the vectors $u_t\in \RR_{+}^{n_t}$, for $t = 1, 2, \ldots, \ccT$, are given by
\begin{equation} \label{eq:u_multi_omt}
u_t = \begin{cases} \exp(  \lambda_t/\epsilon),& \text{ if } t \in \Gamma \\
\ett,& \text{ else,} \end{cases}
\end{equation}
where $\lambda_t\in \RRext^{n_t}$ for $t \in \Gamma$ are optimal dual variables  in the dual problem
of \eqref{eq:omt_multi_regularized}. This dual problem takes the form
\begin{equation} \label{eq:multi_omt_dual}
\maxwrt[\lambda_t\in \RRext^{n_t},\; t\in \Gamma] - \epsilon \langle\bK, \bU \rangle + \sum_{t \in \Gamma} \lambda_t^T \mu_t,
\end{equation}
where $\bU$ depends on $\{\lambda_t\}_{t\in \Gamma}$ as specified in \eqref{eq:U} and \eqref{eq:u_multi_omt}.
For details the reader is referred to, e.g., \cite{elvander19multi, benamou2015bregman}.

The Sinkhorn scheme for finding $\bU$ in \eqref{eq:U} 
is to iteratively update $u_t$ according to
\begin{equation} \label{eq:sinkhorn_multi}
u_t \leftarrow u_t \odot \mu_t ./ P_t(\bK \odot \bU),
\end{equation}
for all $t\in\Gamma$. This scheme may for instance be derived as Bregman projections \cite{benamou2015bregman} or a block coordinate ascend in the dual \eqref{eq:multi_omt_dual}, \cite{karlsson2017generalized,elvander19multi,tseng1990dual}. 
As a result, global convergence of the Sinkhorn scheme \eqref{eq:sinkhorn_multi} is guaranteed \cite{BauLew00, tseng1990dual, luo1992convergence}.
The computational bottleneck of the Sinkhorn iterations \eqref{eq:sinkhorn_multi} is computing the projections $P_t(\bK \odot \bU)$, for $t \in \Gamma$, which in general scales exponentially in $\ccT$. 
In fact, even storing the tensor $\bM$ is a challenge as it  consists of $\prod_{t=1}^\ccT n_t$ elements.
However, in many cases of interest, structures in the cost tensors can be exploited to perform the sum operations in \eqref{eq:proj_discrete} in an appropriate order, which makes the computation of the projections feasible \cite{elvander19multi, benamou2015bregman, haasler20trees, haasler20steering, haasler20pgm}.
More precisely, in many applications the tensor $\bK\odot \bU$ factorizes such that it can be described by a graph $\Gomt=(\Vomt,\Eomt)$, where the vertices $\Vomt$ correspond to the tensor marginals and its dependencies are described by the set of edges $\Eomt$.
The projections \eqref{eq:proj_discrete} can then be computed efficiently by first eliminating the variables, i.e., performing the sum operations, for the vertices that have few dependencies.
For instance, when the tensor $\bK\odot \bU$ factorizes according to a tree structure, the projections \eqref{eq:proj_discrete} can be computed by first eliminating the variables corresponding to the trees leafs and successively moving down the branches. Computing the projections requires then only matrix-vector multiplications, where the matrices are at most of size $\max_t (n_t)$ \cite{haasler20trees, haasler20pgm}.
In the case of more complex graphs a similar approach can be utilized, but computations become more expensive. For instance, in case the graph is a cycle the complexity is increased by a factor of $\max_t (n_t)$ as compared to the tree setting \cite{benamou2015bregman, haasler20steering}.

\section{Network flow problems via optimal transport.} \label{sec:nf_ot}

In this section we introduce a reformulation of the dynamic minimum-cost flow problem as a multi-marginal optimal transport problem \eqref{eq:omt_multi_discrete}.
In the single-commodity case this optimal transport problem has a path-structure. The multi-commodity case can be expressed as several single-commodity problems, which are coupled through the capacity constraints. Alternatively, this can be set up as one multi-marginal optimal transport problem, where the cost function decouples as a graph that contains cycles.

\subsection{The dynamic minimum-cost flow problem.} \label{sec:single_com}

Let $\Nflow_{\rm exp}$ be the time-expansion of the network $\Nflow$ for the time steps $t=0,\dots,\ccT$, and let $\ccP$ denote the set of feasible paths in $\Nflow_{\rm exp}$.
In order to solve an arc-chain formulation of a flow-problem on this network, one has to identify all paths in this set. 
Clearly, the set of feasible paths $\ccP$ is a subset of the set $\tilde \ccP = \{ (i_1,\dots,i_\ccT) : i_t\in \Eflow_t \text{ for } t=1,\dots,\ccT \}$, which contains all combinations of $\ccT$ edges in $\Eflow$.
In fact, the set $\tilde \ccP$ is generally much larger than $\ccP$, since it lifts the set of feasible paths to the set of all ``paths" possible from purely combinatorial considerations (ignoring the graph structure).

\begin{figure}
\centering
	\begin{tikzpicture}
	%\footnotesize
	\tikzstyle{main}=[rectangle, minimum size = 5pt, thick, draw =black!80, node distance = 30pt and 100pt]
	
	\node[main,fill=black!10] (t0n1) {}; 	
	\node[main,fill=black!10] (t0n2) [below =of t0n1] {};  
	\node[main,fill=black!10] (t0n3) [below =of t0n2] {}; 
	\node[draw=none] (t0) [below =10pt of t0n3] {\footnotesize $t=0$}; 
 
	\node[main,fill=black!10] (t1n1) [right =of t0n1] {}; 	
	\node[main,fill=black!10] (t1n2) [below =of t1n1] {};  
	\node[main,fill=black!10] (t1n3) [below =of t1n2] {};
	\node[draw=none] (t1) [below =10pt of t1n3] {\footnotesize $t=1$}; 
		
 	\node[main,fill=black!10] (t2n1) [right =of t1n1] {}; 	
	\node[main,fill=black!10] (t2n2) [below =of t2n1] {};  
	\node[main,fill=black!10] (t2n3) [below =of t2n2] {}; 
	\node[draw=none] (t2) [below =10pt of t2n3] {\footnotesize $t=2$}; 	
	
	\node[main,fill=black!10] (t3n1) [right =of t2n1] {}; 	
	\node[main,fill=black!10] (t3n2) [below =of t3n1] {};  
	\node[main,fill=black!10] (t3n3) [below =of t3n2] {}; 
	\node[draw=none] (t3) [below =10pt of t3n3] {\footnotesize $t=3$}; 	
	
	\draw[ -latex] (t0n1) --  (t1n2);
	\draw[ -latex] (t0n2) -- (t1n3);
	\draw[ -latex] (t0n3) -- (t1n1); 
	\draw[ -latex] (t0n3) -- (t1n2);
	
	\draw[ -latex] (t1n1) --  (t2n2);
	\draw[ -latex] (t1n2) -- (t2n3);
	\draw[ -latex] (t1n3) -- (t2n1); 
	\draw[ -latex] (t1n3) -- (t2n2); 
	
	\draw[ -latex] (t2n1) --  (t3n2);
	\draw[ -latex] (t2n2) -- (t3n3);
	\draw[ -latex] (t2n3) -- (t3n1); 
	\draw[ -latex] (t2n3) -- (t3n2); 

	\draw[dashed, thick] (52.5pt,-37.5pt) ellipse (30pt and 70pt)  ;
	\draw[dashed, thick] (160pt,-37.5pt) ellipse (30pt and 70pt);
	\draw[dashed, thick] (267.5pt,-37.5pt) ellipse (30pt and 70pt);
	
	\node[draw=none] at (52.5 pt, -90pt) (mu1) {$P_1(\bM)$};
	\node[draw=none] at (160 pt, -90pt) (mu1) {$P_2(\bM)$};
	\node[draw=none] at (267.5 pt, -90pt) (mu1) {$P_3(\bM)$};

	\end{tikzpicture}
	\caption{Illustration of the optimal transport tensor $\bM$ in the time-expanded network from Figure~\ref{fig:time_exp}. The tensors marginal $P_t(\bM)$ describes the distribution of flow over the edges in the time-interval $(t-1,t)$.} \label{fig:time_exp_tensor}
\end{figure}
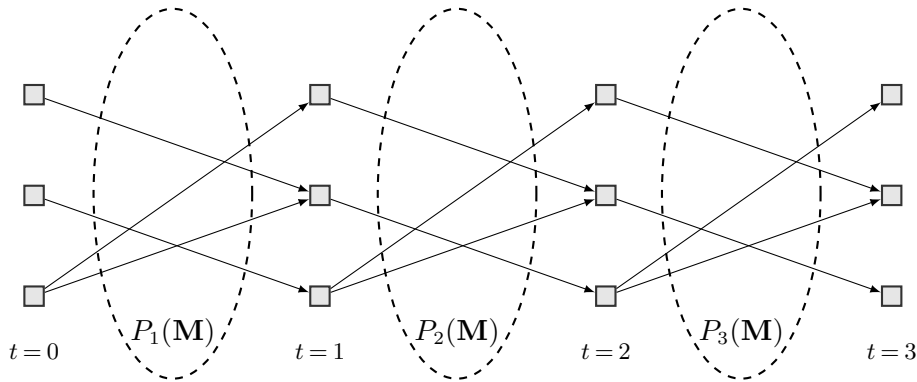

However, using this representation, the network flow can be described by a tensor $\bM\in \RR_{+}^{n^\ccT}$, where $n=|\Eflow|$, and where the element $\bM_{i_1,\dots,i_\ccT}$ denotes the amount of flow on the path $(i_1,\dots,i_\ccT)$. 
The vector $P_t(\bM) \in \RR_{+}^n$, where the projection operator is defined as in \eqref{eq:proj_discrete}, then describes the flow distribution over the edges between time $t-1$ and $t$, as illustrated in Figure~\ref{fig:time_exp_tensor}. That is, its element $P_t(\bM)_i$ denotes the amount of flow over edge $i\in \Eflow_t$.

\begin{figure}
	\centering
	\begin{tikzpicture}
	\tikzstyle{main}=[circle, minimum size = 45pt, thick, draw =black!80, node distance = 40pt]
	\node[main,fill=black!10] (mu0)  {\footnotesize $P_1(\bM)$};  
	\node[main] (mu1) [right=of mu0] {\footnotesize $P_{2}(\bM)$};
	\node[] (mu1T) [right=of mu1] {};  
	\node[main] (muTm1) [right=of mu1T] {\footnotesize $P_{\!\ccT\!-\!1\!}(\bM)\!$};  
	\node[main,fill=black!10] (muT) [right=of muTm1] {\footnotesize $P_{\ccT}(\bM)$};
	
	\draw[->, -latex, thick] (mu0) -- node[below] {$C$} (mu1);
	\draw[->, -latex, thick] (muTm1) -- node[below] {$C$} (muT);
	
	\node[node distance = 1mm] (c1) [below=of mu0] {$c$};
	\node[node distance = 1mm] (c2) [below=of mu1] {$c$};
	\node[node distance = 1mm] (cTm1) [below=of muTm1] {$c$};
	\node[node distance = 1mm] (cT) [below=of muT] {$c$};
	\draw[loosely dotted, very thick] (mu1) -- (muTm1); 
	\end{tikzpicture}
	\caption{Illustration of the path graph for the single-commodity network flow problem. Gray and white circles describe equality and inequality constrained marginals, respectively. As described by \eqref{eq:T_cost}, the costs $c$ are acting on the marginals, and the costs $C$ are acting on the bi-marginals. }
	\label{fig:path_graph}
\end{figure}
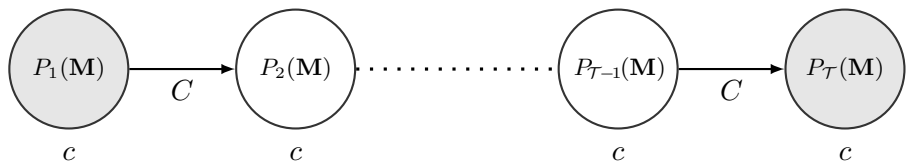

Similarly, the evolution of flow between time intervals $(t_1-1,t_1)$ and $(t_2-1,t_2)$ is described by the bi-marginal projections $P_{t_1,t_2}(\bM) \in \RR_{+}^{n\times n}$, which are defined as
\begin{equation} \label{eq:proj_discrete_bi}
\left( P_{t_1, t_2} (\bM) \right)_{i_{t_1} i_{t_2}} = \sum_{i_{1},\dots, i_\ccT \setminus \{i_{t_1},i_{t_2}\}} \bM_{ i_{1} \dots i_\ccT}.
\end{equation} 
That is, the element $( P_{t_1, t_2} (\bM) )_{ij}$ describes the amount of flow that is in edge $i$ at time $t_1$ and that is in edge $j$ at time $t_2$.
Let $c\in \RR_+^n$, where $c_{i}$ denotes the cost of a unit flow on edge $i\in \ccE$, and let $C \in \RRext_+^{n\times n}$ encode the network topology, i.e., $C_{ij}=0$ if edge $i$ leads to%
\footnote{That is, the second vertex of edge $i$ is the first vertex of edge $j$ in the network $\mathcal{N}$.}
edge $j$, and $C_{ij}=\infty$ otherwise.
Then we define the cost of a transport plan $\bM$ as
\begin{equation} \label{eq:T_cost}
 \sum_{t=1}^{\ccT} c^T P_t(\bM) +  \sum_{t=1}^{\ccT-1} \tr(C^TP_{t,t+1}(\bM))  = \langle \bC, \bM \rangle,
\end{equation}
where the tensor $\bC\in \RRext_+^{n\times n \ldots \times n}$ is defined as
\begin{equation} \label{eq:C_tensor_single_commodity}
\bC_{i_1\dots i_\ccT} = \sum_{t=1}^\ccT c_{i_t} + \sum_{t=1}^{\ccT-1} C_{i_t i_{t+1}}.
\end{equation}
Note that $\langle \bC, \bM \rangle  = \infty $ here means that the transport plan contains paths that are not consistent with the network structure, i.e., that for some $t \in \{1, \ldots, \ccT-1\}$ and some $(i,j) \not \in \ccE$, $(P_{t,t+1}(\bM))_{ij} > 0$.
The structure of the cost function \eqref{eq:T_cost} can be illustrated by the path-graph in Figure~\ref{fig:path_graph}.

Let $\mu_1$ and $\mu_{\ccT}$ be the supply and demand distributions, respectively. That is, $(\mu_1)_i = r_i^+$, for $i\in\ccS^+$, and $0$ otherwise, and $(\mu_\ccT)_i = r^-_i$ for  $i \in \ccS^-$, and  $0$ otherwise. Moreover, let $d\in\RR_{+}^n$ encode the capacity constraints of the network, that is $d_i$ is the flow capacity on edge $i\in \Eflow$. These supply, demand, and capacity constraints can be encoded as equality and inequality constraints on the flow distributions over the edges $P_t(\bM)$.
Based on this, we formulate the linear program
\begin{subequations} \label{eq:singlecommodity}
	\begin{align}
	\min_{\bM\in \RR_+^{ n^\ccT}} \quad &   \langle \bC, \bM \rangle  \label{eq:singlecommodity_obj} \\ 
	\mbox{subject to } \qquad &  P_{1}(\bM)=\mu_1, \label{eq:singlecommodity_s}\\
	&  P_{\ccT}(\bM)= \mu_\ccT  \label{eq:singlecommodity_d} \\
	&P_{t}(\bM)\le d, \qquad \mbox{for } \qquad t=2,\ldots, \ccT-1 \label{eq:singlecommodity_inequality}.
	\end{align}
\end{subequations}
This problem is equivalent to the dynamic minimum-cost network flow problem \eqref{eq:flow_dynamic_single-commodity} in the sense described in the following theorem.

\begin{theorem} \label{thm:equiv_sc} The dynamic minimum-cost network flow problem \eqref{eq:flow_dynamic_single-commodity} and problem \eqref{eq:singlecommodity} correspond to each other in the following sense.
	\begin{enumerate}
\item Assume that \eqref{eq:singlecommodity} has a feasible solution with finite objective value.
Then it has a finite optimal value, and \eqref{eq:flow_dynamic_single-commodity} has the same optimal value.
Moreover, if $\bM$ is an optimal solution of \eqref{eq:singlecommodity}, then there is an optimal solution $\{ x_p : p\in \ccP\}$ of \eqref{eq:flow_dynamic_single-commodity} such that 
\begin{equation}\label{eq:transport_equal}
\bM_{i_1 \dots i_\ccT}=\begin{cases}
x_p & \mbox{ for } (i_1,\dots,i_\ccT) \in \ccP, \mbox{ where } p=(i_1,\dots,i_\ccT)\\
0 & \mbox{ for } (i_1,\dots, i_\ccT) \in \tilde \ccP \setminus \ccP.
\end{cases}
\end{equation}
\item Assume that there is a feasible solution to \eqref{eq:flow_dynamic_single-commodity}. Then it has a finite optimal value, and problem \eqref{eq:singlecommodity} has the same optimal value. Moreover, if $\{ x_p : p\in \ccP\}$ is an optimal solution of \eqref{eq:flow_dynamic_single-commodity}, then there is an optimal solution $\bM$ of \eqref{eq:singlecommodity} such that \eqref{eq:transport_equal} holds.
	\end{enumerate}	
\end{theorem}

\proof{Proof:}
	First, note that the amount of flow on edge $e \in \Eflow$ between time $t-1$ and $t$ is given in the optimal transport formulation \eqref{eq:singlecommodity} by
	\begin{equation} \label{eq:edge_flow_ot}
	\sum_{i \in \tilde P , i_t= e} \bM_{i_1 \dots i_\ccT} = P_t(\bM)_e.
	\end{equation}	
	Thus the flow distribution over $\Eflow_t$ is exactly the projection $P_t(\bM)$ as defined in \eqref{eq:proj_discrete}.
	Then, with
	$$ (\mu_1)_i = \begin{cases} r^+_i, & i \in \ccS^+ \\ 0, & \text{ otherwise,} \end{cases} \qquad
	(\mu_\ccT)_i = \begin{cases} r^-_i, & i \in \ccS^- \\ 0, & \text{ otherwise,} \end{cases} $$
	the set of constraints \eqref{eq:flow_dynamic_single-commodity_s}-\eqref{eq:flow_dynamic_single-commodity_d} and \eqref{eq:singlecommodity_s}-\eqref{eq:singlecommodity_d} both restrict the respective problems to paths that satisfy the supply and demand constraints.
	In the formulation \eqref{eq:flow_dynamic_single-commodity} the total flow on edge $e \in \Eflow_t$ is given by
	\begin{equation} \label{eq:edge_flow}
	\sum_{p \in \ccP} \delta_{e= p_t} x_p,
	\end{equation}	
	and thus the inequality constraints \eqref{eq:flow_dynamic_single-commodity_capacity} and \eqref{eq:singlecommodity_inequality} restrict the flows in the respective problems to the same capacity constraints.
	Moreover, note that $(P_{t,t+1}(\bM))_{ij}$ describes the amount of flow moving from edge $i\in\Eflow_t$ to edge $j\in\Eflow_{t+1}$. Therefore, the objective \eqref{eq:singlecommodity_obj} is finite if and only if $\bM_{i_1\dots i_\ccT}=0$ for all $(i_1,\dots,i_\ccT) \in \tilde \ccP \setminus \ccP$.
	Now, by associating the amount of flow on edge $i \in \ccE_t$ with \eqref{eq:edge_flow_ot} and \eqref{eq:edge_flow}, respectively, the cost of a feasible flow plan, i.e., a plan that satisfies $\bM_{i_1\dots i_\ccT}=0$ for all $(i_1,\dots,i_\ccT) \in \tilde \ccP \setminus \ccP$, can be written in the two formulations as
	\begin{equation}
	\sum_{p \in \ccP} c_p x_p = \sum_{p \in \ccP} \Big( \sum_{t=1}^\ccT \sum_{e\in \Eflow} \delta_{e=p_t} c_e \Big) x_p =   \sum_{e\in \Eflow} \sum_{t=1}^\ccT \Big( \sum_{p \in \ccP}  \delta_{e=p_t} x_p \Big) c_e = \sum_{t=1}^\ccT \sum_{e\in \Eflow} P_t(\bM)_e c_e = \sum_{t=1}^\ccT c^T P_t(\bM).
	\end{equation}
	This completes the proof. \hfill \Halmos
\endproof

Comparing problem \eqref{eq:singlecommodity} to problem \eqref{eq:flow_dynamic_single-commodity}, we have expanded the set of optimization variables by adding a large number of infeasible paths.
However, the novel formulation  \eqref{eq:singlecommodity} is structured as a multi-marginal optimal transport problem as in \eqref{eq:omt_multi_discrete}, which opens up for efficiently computing an approximate solution. 
In particular, the structure of problem \eqref{eq:singlecommodity} can be described by the path graph in Figure~\ref{fig:path_graph}.
Although problem \eqref{eq:singlecommodity} lifts the set of optimization variables in \eqref{eq:flow_dynamic_single-commodity} from the set of feasible paths to the set of all combinatorially possible paths in the network, the infinite values in the tensor \eqref{eq:C_tensor_single_commodity} restrict the problem to the set of feasible paths as in \eqref{eq:flow_dynamic_single-commodity}.
\begin{remark}
	The second term in \eqref{eq:T_cost} is needed only to restrict the solution of problem \eqref{eq:singlecommodity} to the set of feasible paths $\ccP$. Naturally this could instead be imposed as a set of hard constraints $P _{t,t+1} (\bM) \leq E$, for $t=1,\dots, \ccT-1$, where $E_{ij}=\infty$ if edge $i$ leads to edge $j$, and $E_{ij}=0$ otherwise. Instead, we choose to use the penalty terms in \eqref{eq:T_cost} for computational reasons. In Section~\ref{sec:OTgraph} we develop a scheme, which is based on the methods introduced in Section~\ref{sec:omt}, i.e., solving the dual of a regularization of problem \eqref{eq:singlecommodity}. Note that adding more hard constraints to \eqref{eq:singlecommodity} leads to a larger number of dual variables, which makes it more expensive to solve the regularized dual problem. We thus impose the network structure through the penalty terms in \eqref{eq:T_cost}, which yields a dual problem with considerably fewer variables.
Moreover, infinite values in $C$	induce sparsity to the tensor $\bK= \exp(-\bC/\epsilon)$ in \eqref{eq:MKU}, which can be exploited when computing the projections \eqref{eq:proj_discrete} needed for the Sinkhorn scheme.
\end{remark}

\subsection{The dynamic multi-commodity minimum-cost flow problem.} \label{sec:mc_form}

In this section we extend the optimal transport formulation of the dynamic minimum-cost network flow problem from Section~\ref{sec:single_com} to the multi-commodity setting.

Assume that there are $L$ different commodities present in the network $\Nflow$, and each of these is assigned an initial distribution $\mu_1^\ell$ and a final distribution $\mu_\ccT^\ell$, for $\ell=1,\dots,L$.
For each commodity we define a cost vector $c_\ell \in \RR^n$, where $(c_\ell)_i$ denotes the cost of a unit flow of commodity $\ell$ on edge $i \in \Eflow$.
As in the single-commodity case in Section~\ref{sec:single_com}, the network structure is imposed by a matrix $C\in \RRext_{+}^{n\times n}$, and the total flow capacity is bounded on all edges, and described by a vector $d\in \RR_{+}^n$.
One way to formulate an optimal transport problem for the multi-commodity flow is to describe each commodities flow by a mass transport tensor $\bM^\ell$, for $\ell=1,\dots,L$.
Then each of these transport tensors has to satisfy the respective supply and demand constraints \eqref{eq:singlecommodity_s}-\eqref{eq:singlecommodity_d}, and its cost is given by $\langle \bC^\ell, \bM^\ell \rangle$ 
as defined in \eqref{eq:T_cost}.
The capacity constraints in the network need to hold for the sum of all commodity flows, i.e., the sum of the projections $P_t( \bM^\ell)$ over all commodities $\ell=1,\dots,L$.
The dynamic multi-commodity minimum-cost flow problem \eqref{eq:flow_dynamic_multi-commodity} can therefore be written as
\begin{subequations} \label{eq:multicommodity_tensor_problem}  \noeqref{eq:multicommodity_tensor_problem_s,eq:multicommodity_tensor_problem_d}
	\begin{align} 
	\minwrt[ \bM^1,\dots,\bM^L \in \RR_{+}^{n^{\ccT}}  ] \ &  \sum_{\ell=1}^L  \langle \bC^\ell, \bM^\ell \rangle %\sout{ \quad  T_{c_\ell,C}(\bM^\ell) }
	 \label{eq:multicommodity_tensor_problem_obj} \\
	\text{subject to } \ & \ P_{1} ( \bM^\ell ) = \mu_1^\ell, \quad \text{ for } \ell=1,\dots,L, \label{eq:multicommodity_tensor_problem_s} \\
	& \ P_{\ccT} ( \bM^\ell ) = \mu_\ccT^\ell, \quad \text{ for } \ell=1,\dots,L, \label{eq:multicommodity_tensor_problem_d} \\
	& \ \sum_{\ell=1}^L P_t( \bM^\ell) \leq d, \quad \text{ for }  t=2,\dots,\ccT-1. \label{eq:multicommodity_tensor_problem_cap}
	\end{align}
\end{subequations}
Note here that the $L$ optimal transport problems are each of the form in \eqref{eq:singlecommodity}, and are coupled only through the capacity constraint \eqref{eq:multicommodity_tensor_problem_cap}.

\begin{figure}
	\centering
	\begin{tikzpicture}
	\footnotesize
	%\scriptsize
	\tikzstyle{main}=[circle, minimum size = 45pt, thick, draw =black!80, node distance = 40pt]
	\tikzstyle{obs}=[circle, minimum size = 45pt, thick, draw =black!80, node distance = 50pt and 30pt ]
	%\node[main,fill=black!10] (mu0) {$\mu_0$};
	\node[main,fill=black!10] (mu0) {$P_0(\bM)$};
	
	\node[] (phi1c) [below=of mu0] {};  
	\node[obs] (mu2) [left=of phi1c] {$P_{2}(\bM)$};  
	\node[obs,fill=black!10] (mu1) [left=of mu2] {$P_{1}(\bM)$};  
	\node[obs] (muS1) [right=of phi1c] {$\!P_{\!\ccT\!-\!1\!}(\bM)\!$};
	\node[obs,fill=black!10] (muS) [right=of muS1] {$P_{\ccT}(\bM)$};
	
	\draw[->, -latex, thick] (mu0) -- node[above left] {$P_{0,1}(\bM)=R^{0,1}$} (mu1);
	\draw[->, -latex, thick] (mu0) -- node[left] {$\CL$} (mu2);
	\draw[->, -latex, thick] (mu0) -- node[left] {$\CL$} (muS1);
	\draw[->, -latex, thick] (mu0) -- node[above right] {$P_{0,\ccT}(\bM)=R^{0,\ccT}$} (muS);
	\draw[->, -latex, thick] (mu1) -- node[below] {$C$} (mu2);
	\draw[->, -latex, thick] (muS1) -- node[below] {$C$} (muS);
	
	\draw[loosely dotted, very thick] (mu2) -- (muS1); 
	\end{tikzpicture}
	\caption{Illustration of the dynamic multi-commodity minimum cost flow problem \eqref{eq:multicommodity}. Gray and white circles describe equality and inequality constrained marginals, respectively. } 
	\label{fig:multicommodity}
\end{figure}
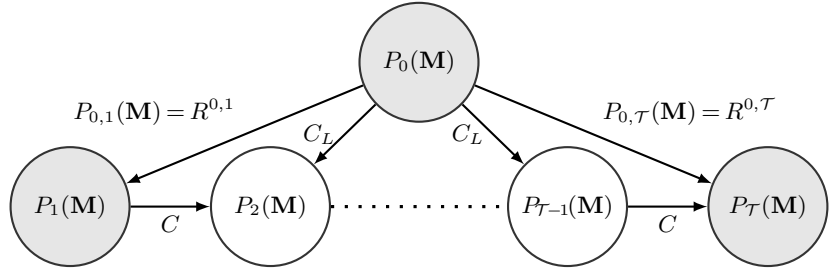

We will now bring problem \eqref{eq:multicommodity_tensor_problem} on a form similar to a multi-marginal optimal transport problem \eqref{eq:omt_multi_discrete}, i.e., a formulation containing only one mass transport tensor.
This is done by combining all information from the $\ccT$-mode transport plans $\bM^\ell \in \RR_{+}^{n \times \dots \times n}$, for $\ell=1,\dots,L$, to a new mass transport tensor $\bM \in \RR_{+}^{L\times n \times \dots \times n}$ with $\ccT+1$ modes.
That is, we let its element $\bM_{\ell, i_{1} \dots i_\ccT}$ describe the amount of flow of commodity $\ell$ over the path $i_1,\dots,i_\ccT$.
Accordingly, for the added mode in the tensor we introduce a marginal $\mu_0\in \RR_+^L$, where $(\mu_0)_\ell=\ett^T \mu_1^\ell =\ett^T \mu_\ccT^\ell$ denotes the total supply and demand of commodity $\ell\in L$.
The initial and final distributions for the commodities can then be summarized in two matrices $R^{(0,1)},R^{(0,\ccT)} \in \RR_+^{L\times n}$, defined as $R^{(0,1)}=(\mu_1^1, \mu_1^2,\ldots, \mu_1^L)^T $ and $R^{(0,\ccT)}=(\mu_\ccT^1, \mu_\ccT^2,\ldots, \mu_\ccT^L)^T$.
In particular, with this construction it holds that $R^{(0,1)} \ett = R^{(0,\ccT)} \ett = \mu_0$.
Moreover, define a matrix $\CL \in \RR^{L\times n}$ as $\CL=(c_1,c_2,\dots,c_L)^T$, that is $(\CL)_{\ell,i}$ denotes the cost for commodity $\ell\in L$ to be on edge $i \in \ccE$.
This setup is illustrated in Figure~\ref{fig:multicommodity}. 

Note that the objective function \eqref{eq:multicommodity_tensor_problem_obj} can be written as
\begin{equation}
\sum_{t=2}^{\ccT-1} \tr(\CL^TP_{0,t}(\bM))+  \sum_{t=1}^{\ccT-1} \tr(C^TP_{t,t+1}(\bM)) =  \langle \bC, \bM \rangle,
\end{equation}
where the cost tensor $\bC\in \RRext_+^{L \times n^\ccT}$ is given by
\begin{equation} \label{eq:cost_multicommodity}
\bC_{i_0\dots i_\ccT} = \sum_{t=2}^{\ccT-1} (\CL)_{i_0 i_t} +  \sum_{t=1}^{\ccT-1} C_{i_t i_{t+1}}.
\end{equation}
Thus, the dynamic multi-commodity minimum-cost network flow problem \eqref{eq:multicommodity_tensor_problem} can be expressed as
\begin{equation}  \label{eq:multicommodity}
\begin{aligned}
\minwrt[\bM\in \RR_+^{L\times n^\ccT}] \ & \  \langle \bC, \bM \rangle
 \\ 
\mbox{subject to } \
&\ P_{0,1}(\bM)=R^{(0,1)},  \\
&\ P_{0,\ccT}(\bM)=R^{(0,\ccT)}, \\
&\ P_{t}(\bM)\le d_t, \qquad \mbox{for } \qquad t=2,\ldots, \ccT-1
\end{aligned}
\end{equation}

Utilizing the result in Theorem~\ref{thm:equiv_sc} we have now proved that the solution to \eqref{eq:multicommodity} and the dynamic multi-commodity minimum-cost network flow problem \eqref{eq:flow_dynamic_multi-commodity} are equivalent, as summarized in the following Theorem.
\begin{theorem} \label{thm:equiv_mc}
The dynamic minimum-cost network flow problem \eqref{eq:flow_dynamic_multi-commodity} and problem \eqref{eq:multicommodity} correspond to each other in the following sense.
	\begin{enumerate}
\item Assume that \eqref{eq:multicommodity} has a feasible solution with finite objective value.
Then \eqref{eq:multicommodity} has a finite optimal value, and \eqref{eq:flow_dynamic_multi-commodity} has the same optimal value.
Moreover, if $ \bM$ is an optimal solution of \eqref{eq:multicommodity}, then there is an optimal solution $\{ x^\ell_p : p\in \ccP^\ell, \ell=1,\dots, L\}$ of \eqref{eq:flow_dynamic_multi-commodity} such that 
\begin{equation}\label{eq:transport_equal_multi}
\bM_{\ell, i_1 \dots i_\ccT}=\begin{cases}
x^\ell_p & \mbox{ for } (i_1,\dots,i_\ccT) \in \ccP^\ell, \mbox{ where } p=(i_1,\dots,i_\ccT)\\
0 & \mbox{ for } (i_1,\dots, i_\ccT) \in \tilde \ccP \setminus \ccP^\ell.
\end{cases}
\end{equation}
\item Assume that there is a finite feasible solution to \eqref{eq:flow_dynamic_multi-commodity}. Then it has a finite optimal value, and problem \eqref{eq:multicommodity} has the same optimal value. Moreover, if $\{ x^\ell_p : p\in \ccP^\ell, \ell= 1, \dots, L \}$ is an optimal solution of \eqref{eq:flow_dynamic_multi-commodity}, then there is an optimal solution $\bM$ of \eqref{eq:multicommodity} such that \eqref{eq:transport_equal_multi} holds.
	\end{enumerate}	
\end{theorem}

\subsection{Generalizations.}
\label{subsec:discussion}

In this section we have introduced novel formulations for dynamic minimum-cost network flow problems based on the optimal transport framework.
We will now discuss a few modifications and generalizations of the proposed problems \eqref{eq:singlecommodity} and \eqref{eq:multicommodity}, and show that the proposed formulations in fact provide a highly flexible framework for dynamic network flow problems.

An advantage of our framework is that, due to the fact that the network structure is imposed by the cost matrix $C$, a time-varying network can be modelled in a straightforward way. Namely, the matrix $C$ can simply be replaced by a set of time-dependent matrices $C_t$, for $t=1,\dots,\ccT-1$, where $C_t$ encodes the network topology in the interval $(t,t+1)$.
Moreover, based on the formulation \eqref{eq:multicommodity_tensor_problem}, where each commodity is described by a separate transport tensor, one can extend the problem to the setting, where different commodities enter and leave the network at different times. 
In fact, the computational methods derived in this work can easily be modified to this setting, as we will argue in Remark~\ref{rem:multi_tensor_scheme}.

 In some applications, for instance in traffic flow problems, where edges and nodes describe streets and junctions, respectively, it is natural to allow for intermediate storage on the edges.
This can be easily incorporated in our framework by letting $C_{ii}$ denote the cost for staying on edge $i \in \Eflow$.
It should be noted that in this case the cost $c_i^\ell$ denotes the cost for commodity $\ell$ to use edge $i \in \Eflow$, and not the cost for traveling between the two vertices. That is, the cost accumulates if flow remains on an edge for several time intervals, which is useful, e.g., in traffic routing problems, where the cost model should take the travel time of agents into account.
However, we can achieve a cost that does not accumulate in the case where all commodities are described by the same cost $c_i= c_i^\ell$ for all $\ell=1,\dots,L$, by defining a negative cost $C_{ii} = - c_i$ for staying on the edge $i\in\Eflow$.

A more classical setting in network flow problems is to allow for storage in the vertices.
One way to include this in the presented framework is to augment the support of the modes of the mass transport tensor by the set of vertices, i.e., by letting $n= |\Eflow| + |\Vflow|$.
In particular, in the multi-commodity problem \eqref{eq:multicommodity} the mass transport tensor is then of the size $\bM \in \RR_{+}^{L \times (|\Eflow| + |\Vflow|)^\ccT }$, and the distributions are of the size $\mu_t \in \RR_{+}^{|\Eflow| + |\Vflow|}$, for $t=1,\dots,\ccT$.
Analogously to before, the network structure is imposed by the cost matrices $C \in \RR_{+}^{(|\Eflow| + |\Vflow|) \times (|\Eflow| + |\Vflow|)}$, i.e., we define $C_{ij}= 0$ if $i \in \{\Eflow \cup \Vflow \}$ %leads%
%\footnote{A vertex leads to all edges it is the first element of, and to itself.}
is adjacent%
\footnote{A vertex is adjacent to all edges it connects to, and to itself.} 
to $j \in \{\Eflow \cup \Vflow \}$, and $C_{ij}= \infty$ otherwise.
Similarly, the definition of the cost $\CL$ and the capacities $d$ can be extended to the vertices.
It is worth noting that this extension of the state space also allows for defining the set of sinks and sources on the vertices instead of the edges.

Another extension of the formulation, of particular interest for traffic routing problems, is the setting where the sinks and sources are defined on nodes, but intermediate storage is only allowed in the sinks and sources, and agents are not permitted to enter sources, or leave sinks. In this case, we let $n=|\Eflow|+|\ccS^+| + |\ccS^-|$, and define the network structure through the cost matrix as follows
\begin{equation}\label{eq:cost_traffic_routing}
C_{ij} = \begin{cases} 0, & \text{if } i\in \{\Eflow \cup \ccS^+\} \text{ is adjacent to } j \in \{ \Eflow \cup \ccS^+  \cup \ccS^- \} \\
0, & \text{if } i\in \{\Eflow \cup \ccS^+ \cup \ccS^-\} \text{ is adjacent to } j \in \{ \Eflow \cup \ccS^-   \} \\
\infty, & \text{otherwise.}
\end{cases}
\end{equation}
%\end{remark}

A final extension worth mentioning is the possibility of introducing commodity-dependent capacity constraints \cite{kennington1978survey, gendron1999multicommodity}. This may be done by introducing the set of constraints $P_{0,t}(\bM) \leq D^{(0,t)}$ for $t=2,\dots, \ccT-1$, with capacity matrices $D^{(0,t)} \in \RRext_{+}^{L \times n}$, where $D^{(0,t)}_{\ell i}$ denotes the capacity of commodity $\ell$ on edge $i\in \Eflow$.

\section{The graph-structured multi-marginal optimal transport problem.} \label{sec:OTgraph}
In this section we define the general graph structured optimal transport problem and develop methods to solve the corresponding entropy regularized problem, see also \cite{haasler20trees, haasler20pgm, benamou2015bregman, altschuler2020polynomial}. We will also consider the dynamic flow problems in Section~\ref{sec:nf_ot} in detail, and show how to exploit the graph-structures in order to derive efficient methods.

We have noted that the network flow problems \eqref{eq:singlecommodity} and \eqref{eq:multicommodity} can be seen as multi-marginal optimal transport problems with the underlying graph-structures in Figure~\ref{fig:path_graph} and Figure~\ref{fig:multicommodity}. 
In particular, we let each mode of the transport tensor $\bM$ be associated with a vertex, and let interaction terms be described by edges.
This defines a graph $\Gomt=(\Vomt,\Eomt)$ with vertices $\Vomt$ and edges $\Eomt$.
The interaction terms defining the edges are given by bi-marginal constraints, as in \eqref{eq:multicommodity}, or by bi-marginal cost terms in the cost tensor, i.e., $\bC\in \RR^{n_1\times \cdots \times n_\ccT}$ with
\begin{equation}
´\bC_{i_1 \ldots i_\ccT}=\sum_{(t_1, t_2)\in E}C^{(t_1, t_2)}_{i_{t_1} i_{t_2}}.
\end{equation}
We denote the set of marginals that are constrained by equality and inequality constraints by $\tilde \Vomt_= \subset \Vomt$ and $\tilde \Vomt_\leq \subset \Vomt$, respectively.
Moreover, the set of tuples that are associated with a bi-marginal constraint is denoted by $\tilde \Eomt$.
Thus, the dynamic network flow problems \eqref{eq:singlecommodity} and \eqref{eq:multicommodity} are special cases of the graph-structured optimal transport problem
\begin{equation} \label{eq:omt_multi_graph}
\begin{aligned}
\minwrt[\bM \in \RR_+^{n_1\times\dots\times n_\ccT}] \ & \  \langle \bC, \bM \rangle   \\
\text{ subject to } \ & \ P_t (\bM) = \mu_t,  \qquad \qquad \;\text{ for } t\in \tilde \Vomt_=,\\
& \ P_t (\bM) \leq d_t,  \qquad \qquad \; \text{ for } t\in \tilde \Vomt_\leq, \\
& \  P_{t_1,t_2} (\bM) = R^{(t_1,t_2)},\quad  \text { for } (t_1,t_2) \in \tilde \Eomt,
\end{aligned}
\end{equation}
where $\mu_t,d_t \in \RR_{+}^{n_t}$, and $ R^{(t_1,t_2)} \in \RR_{+}^{n_{t_1} \times n_{t_2}}$.
Following the approach presented in Section~\ref{subsec:sinkhorn} we develop a scheme for approximately solving optimal transport problems of this form.
It is worth noting that the results in Theorem~\ref{thm:solution}, Theorem~ \ref{thm:solution_extension} and Proposition~\ref{prp:sinkhorn} are only based on the structure of the constraints in \eqref{eq:omt_multi_graph}, and thus hold for arbitrary cost tensors $\bC$.
However, to derive the efficient schemes presented in Section~\ref{subsec:scheme_singlecommodity} and Section~\ref{subsec:scheme_singlecommodity} the graph-structures in the objective function have to be exploited.

\subsection{Sinkhorn's method}
In order to apply the approach in Section~\ref{subsec:sinkhorn} we regularize \eqref{eq:omt_multi_graph} with an entropy term \eqref{eq:entropy_term}, which yields the regularized problem
\begin{equation} \label{eq:omt_multi_graph_reg}
\begin{aligned}
\minwrt[\bM \in \RR_+^{n_1\times\dots\times n_\ccT}] \ & \  \langle \bC, \bM \rangle  + \epsilon D(\bM) \\
\text{ subject to } \ & \ P_t (\bM) = \mu_t,  \qquad \qquad \;\text{ for } t\in \tilde \Vomt_=,\\
& \ P_t (\bM) \leq d_t,  \qquad \qquad \;\text{ for } t\in \tilde \Vomt_\leq, \\
& \  P_{t_1,t_2} (\bM) = R^{(t_1,t_2)},\quad  \text { for } (t_1,t_2) \in \tilde \Eomt.
\end{aligned}
\end{equation}
Similarly to the standard multi-marginal optimal transport problem, the solution to \eqref{eq:omt_multi_graph_reg} can be expressed in terms of its optimal dual variables, as the following theorem describes.
\begin{theorem} \label{thm:solution}
	Assume $\bC$ is finite, and the prescribed marginals $\mu_t$ for $t\in \tilde \Vomt_=$, $d_t$ for $t\in \tilde \Vomt_\leq$, and $R^{(t_1,t_2)}$ for $(t_1,t_2) \in \tilde \Eomt$ are strictly positive.
	 Moreover, assume that \eqref{eq:omt_multi_graph_reg} has a feasible solution. 
	Let $\tilde \Vomt = \tilde \Vomt_= \cup \tilde \Vomt_\leq$.
	Then the optimal solution to \eqref{eq:omt_multi_graph_reg} has the structure $\bM=\bK\odot \bU$ where $\bK=\exp(-\bC/\epsilon)$ and 
	\begin{equation} \label{eq:U_tensor}
	\bU_{ i_{1} \ldots  i_{\ccT}}=\left(\prod_{t\in \tilde \Vomt}(u_t)_{i_t}\right) \left(\prod_{(t_1,t_2)\in \tilde \Eomt}U_{i_{t_1}  i_{t_2}}^{(t_1, t_2)}\right),
	\end{equation}
	where $u_t \in \RR_{+}^{n_{t}}$, for $t\in \tilde \Vomt$, and $U^{(t_1,t_2)}  \in \RR_{+}^{n_{t_1} \times n_{t_2}}$, for $(t_1,t_2)\in \tilde \Eomt$.
	
	In particular, $u_t=\exp(-\lambda_{t}/\epsilon)$ and $U^{(t_1,t_2)}= \exp(-\Lambda^{(t_1,t_2)}/ \epsilon)$, where $\lambda_t \in \RR^{n_t}$
	and $\Lambda^{(t_1,t_2)}  \in \RR^{n_{t_1} \times n_{t_2}}$, for $t\in \tilde \Vomt$ and $(t_1,t_2)\in \tilde \Eomt$, respectively,
	are optimal variables for the dual problem of  \eqref{eq:omt_multi_graph_reg}, which is given by
	\begin{equation} \label{eq:omt_multi_graph_dual} 
	\maxwrt[ \substack{ \Lambda^{(t_1,t_2)} \in \RR^{n_{t_1} \times n_{t_2}},\ (t_1,t_2)\in \tilde \Eomt,\\ \lambda_t \in\RR^{n_t},\ t \in \tilde\Vomt_= \\ \lambda_t \in\RR^{n_t}_+,\ t \in \tilde\Vomt_\leq} ] \ \ - \epsilon \langle \bK,  \bU \rangle - \sum_{(t_1,t_2)\in \tilde \Eomt} \langle \Lambda^{(t_1,t_2)}, R^{(t_1,t_2)}  \rangle - \sum_{t \in \tilde \Vomt} \langle \lambda_t, \mu_t \rangle.
	\end{equation}
\end{theorem}
\proof{Proof:}
	Define Lagrange multipliers $\Lambda^{(t_1,t_2)} \in \RR^{n_{t_1} \times n_{t_2}}$, for $(t_1,t_2) \in \tilde \Eomt$, and $\lambda_t \in \RR^{n_t}$, for $t\in \tilde \Vomt$. Moreover, let $\lambda := (\lambda_t)_{t\in \tilde \Vomt}$ and $\Lambda := (\Lambda^{(t_1,t_2)})_{(t_1,t_2) \in \tilde \Eomt}$.
	With these, a Lagrangian of \eqref{eq:omt_multi_graph_reg} is 
	\begin{equation} \label{eq:Lagrangian}
	\mathcal{L}(\bM, \lambda, \Lambda) :=  \langle \bC,  \bM \rangle + \epsilon D(\bM)  + \sum_{(t_1,t_2)\in \tilde \Eomt} \langle \Lambda, P_{t_1,t_2}(\bM) -R^{(t_1,t_2)} \rangle + \sum_{t \in \tilde \Vomt} \langle \lambda_t,  P_t(\bM) - \mu_t \rangle .
	\end{equation}
	The minimum of \eqref{eq:Lagrangian} with respect to $\bM_{i_1\dots i_\ccT}$ is achieved when its derivative vanishes, i.e., when
	\begin{equation}
	\bC_{i_1\dots i_\ccT} + \epsilon \log \left( \bM_{i_1\dots i_\ccT} \right) + \sum_{(t_1,t_2)\in \tilde \Eomt}  \Lambda^{(t_1,t_2)}_{i_{t_1} i_{t_2}} + \sum_{t \in \tilde \Vomt} (\lambda_t)_{i_t}  = 0. 
	\end{equation}
	Thus, the optimal transport tensor is of the form $\bM = \bK \odot \bU$ with $\bK$ and $\bU$ as defined in the theorem.
	Note that the entropy term $D( \bK \odot \bU)$ reads
	\begin{equation}
	\begin{aligned}
	& \sum_{i_1,\dots,i_\ccT} \! \Bigg( \bK_{i_1\dots i_\ccT} \bU_{i_1\dots i_\ccT} \frac{1}{\epsilon} \Big( - \bC_{i_1\dots i_\ccT} -  \sum_{(t_1,t_2)\in \tilde \Eomt}  \Lambda^{(t_1,t_2)}_{i_{t_1} i_{t_2}} - \sum_{t \in \tilde \Vomt} (\lambda_t)_{i_t} \Big)  - \bK_{i_1\dots i_\ccT} \bU_{i_1\dots i_\ccT} + 1 \Bigg) \\
	&= - \frac{1}{\epsilon} \langle \bK \odot \bU, \bC \rangle - \frac{1}{\epsilon} \sum_{(t_1,t_2)\in \tilde \Eomt}  \langle \Lambda^{(t_1,t_2)}, P_{t_1,t_2}(\bK \odot \bU) \rangle - \frac{1}{\epsilon} \sum_{t \in \tilde \Vomt} \langle \lambda_t, P_t(\bK \odot \bU) \rangle - \langle \bK, \bU \rangle + \prod_{t=1}^\ccT n_{t}.
	\end{aligned} 
	\end{equation}
	Thus, plugging $\bM = \bK \odot \bU$ into $\mathcal{L}(\bM, \lambda, \Lambda)$ in \eqref{eq:Lagrangian}, and removing constants, yields
	\begin{equation} \label{eq:dual_obj}
	- \epsilon \langle \bK, \bU \rangle - \sum_{(t_1,t_2)\in \tilde \Eomt} \langle \Lambda^{(t_1,t_2)}, R^{(t_1,t_2)}  \rangle - \sum_{t \in \tilde \Vomt} \langle \lambda_t, \mu_t \rangle.
	\end{equation}
	The dual to \eqref{eq:omt_multi_graph_reg} is to maximize \eqref{eq:dual_obj} with respect to $\Lambda^{(t_1,t_2)}$ for $(t_1,t_2) \in \tilde \Eomt$, and $\lambda_t$ for $t\in \tilde \Vomt$. Finally, given the assumptions, strong duality holds between the primal and the dual problem, see, e.g., \cite[p. 226]{boyd2004convex}. \hfill \Halmos
\endproof

The assumptions in Theorem~\ref{thm:solution} are typically not satisfied for the network flow problems \eqref{eq:singlecommodity} and \eqref{eq:multicommodity}. If the underlying network is not a complete graph, the cost tensor has infinite entries. Moreover, in most flow problems, the sources and sinks are a strict subset of the set of edges, which is modeled by zero entries in the prescribed marginals $\mu_t$, or $R^{(t_1,t_2)}$. The following theorem extends Theorem~\ref{thm:solution} to these cases.
\begin{theorem} \label{thm:solution_extension}
Let $\bC\in \RRext_+^{n_1\times\dots\times n_\ccT}$ and assume that there is a feasible solution $\bM$ of \eqref{eq:omt_multi_graph_reg} for which $\bM_{i_1 \ldots i_\ccT}>0$ if and only if $\bC_{i_1 \ldots i_\ccT}<\infty$, $(\mu_t)_{i_t}>0$, $(d_t)_{i_t}>0$, and $R^{(t_1, t_2)}_{i_{t_1} i_{t_2}}>0$. Then the optimal solution to \eqref{eq:omt_multi_graph_reg} has the structure $\bM=\bK\odot \bU$ where $\bK=\exp(-\bC/\epsilon)$ and $\bU$ factorizes as in \eqref{eq:U_tensor}. 
\end{theorem}

\proof{Proof:}
Define the set of tuples
\begin{equation}
I = \{ (i_1,\dots,i_\ccT) | i_t \in \{1,\dots,n\}, \bC_{i_1 \ldots i_\ccT}<\infty, (\mu_t)_{i_t}>0, (d_t)_{i_t}>0, R^{(t_1, t_2)}_{i_{t_1} i_{t_2}}>0 \}.
\end{equation}
For $(i_1,\dots,i_\ccT) \in I$ we define $\hat C_{i_1 \dots i_\ccT} = C_{i_1 \dots i_\ccT}$, $(\hat \mu_t)_{i_t} = (\mu_t)_{i_t}$, $(\hat d_t)_{i_t}= (d_t)_{i_t}$, and $  \hat R^{(t_1, t_2)}_{i_{t_1} i_{t_2}} =  R^{(t_1, t_2)}_{i_{t_1} i_{t_2}}$.
Consider the problem
\begin{equation} \label{eq:omt_multi_graph_hat}
\begin{aligned}
\minwrt[\hat \bM_{i_1 \dots i_\ccT}, (i_1,\dots,i_\ccT) \in I] \ & \  \sum_{(i_1,\dots,i_\ccT) \in I} \tilde \bC_{i_1 \dots i_\ccT} \hat \bM_{i_1 \dots i_\ccT} + \epsilon D( \hat \bM) \\
\text{ subject to } \ & \ P_t ( \hat \bM) = \tilde \mu_t,  \qquad \qquad \; \text{ for } t\in \tilde \Vomt_=,\\
& \ P_t ( \hat \bM) \leq \hat d_t,  \qquad \qquad \; \text{ for } t\in \tilde \Vomt_\leq, \\
& \  P_{t_1,t_2} ( \hat \bM) = \hat R^{(t_1, t_2)},\quad  \text { for } (t_1,t_2) \in \tilde \Eomt,
\end{aligned}
\end{equation}
where the definition of $D(\bM)$, $P_t(\bM)$ and $P_{t_1,t_2} (\bM)$ is relaxed to the case where the argument is not a tensor.
The proof of Theorem~\ref{thm:solution} can be mirrored for the case where the variable is not a tensor. Thus, the optimal solution to \eqref{eq:omt_multi_graph_hat} can be written as $\hat \bM_{i_1 \ldots i_\ccT} = \hat \bK_{i_1 \ldots i_\ccT}  \hat \bU_{i_1 \ldots i_\ccT}$, where $\hat \bK_{i_1 \ldots i_\ccT} = \exp(- \hat \bC_{i_1 \ldots i_\ccT} / \epsilon)$, and
\begin{equation} \label{eq:U_tensor_hat}
\hat \bU_{i_{1} \ldots i_{\ccT}}=\left(\prod_{t\in \tilde \Vomt}  (\hat u_t)_{i_t}\right) \left(\prod_{(t_1,t_2)\in \tilde \Eomt} \hat U^{(t_1, t_2)}_{i_{t_1} i_{t_2}}\right),
\end{equation}
where $(i_{1},\dots,i_\ccT) \in I$.
Now, define the tensors $\bK =\exp(-\bC/\epsilon)$ and $\bU \in \RR_{+}^{n_1 \times\dots\times n_\ccT}$, which is constructed as in \eqref{eq:U_tensor}, where
\begin{equation} \label{eq:u_augment}
(u_t)_{i_t} = \begin{cases}  (\hat u_t)_{i_t}, & \text{ if it is defined,} \\
0, & \text{ otherwise,} \end{cases} \qquad 
 U^{(t_1, t_2)}_{i_{t_1} i_{t_2}} = \begin{cases}   \hat U^{(t_1, t_2)}_{i_{t_1} i_{t_2}}, & \text{ if it is defined,} \\
0, & \text{ otherwise.} \end{cases}
\end{equation}
Then by construction $\bM = \bK \odot \bU$ is an optimal solution to \eqref{eq:omt_multi_graph_reg}.  \hfill \Halmos
\endproof

The Sinkhorn iterations for problem \eqref{eq:omt_multi_graph_reg} can be derived as a block-coordinate ascend method in the dual problem \eqref{eq:omt_multi_graph_dual}, as summarized in the following proposition.
\begin{proposition} \label{prp:sinkhorn}
Assume \eqref{eq:omt_multi_graph_reg} has a feasible solution as in the assumptions of Theorem~\ref{thm:solution_extension}. Let $\bK=\exp(-\bC/\epsilon)$ and $\bU$ as defined in \eqref{eq:U_tensor}.
Then, the iterative scheme
\begin{subequations} \label{eq:sinkhorn_graph}
\begin{align}
U^{(t_1,t_2)} & \leftarrow U^{(t_1,t_2)} \odot R^{(t_1,t_2)} ./ P_{t_1,t_2} (\bK \odot \bU),  \quad \text{ for } (t_1,t_2)\in \tilde \Eomt  \label{eq:sinkhorn_graph_E} \\
u_t & \leftarrow u_t \odot \mu_t ./ P_t(\bK \odot \bU), \hspace{70pt} \text{ for } t \in \tilde \Vomt_= \label{eq:sinkhorn_graph_Veq} \\
u_t & \leftarrow \min \left( u_t \odot d_t ./ P_t(\bK \odot \bU) \ , \ \ett \right),  \qquad \text{ for } t \in \tilde \Vomt_\leq, \label{eq:sinkhorn_graph_Vineq}
\end{align}	
\end{subequations}
 converges linearly, and in the limit point the optimal solution of \eqref{eq:omt_multi_graph_reg} is given by $\bM = \bK \odot \bU$.
\end{proposition}
\proof{Proof:}
	We first assume that the stronger assumptions from Theorem~\ref{thm:solution} hold. The scheme is derived as a block coordinate ascent method in the dual \eqref{eq:omt_multi_graph_dual}. This is to maximize the objective with respect to one set of dual variables while keeping the other dual variables fixed, i.e., to perform the updates
	\begin{subequations} \label{eq:bca}
		\begin{align}
		\Lambda^{(t_1,t_2)} &\leftarrow \argmaxwrt[\Lambda^{(t_1,t_2)} \in \RR^{n_{t_1} \times n_{t_2}}]   - \epsilon \langle \bK,  \bU \rangle - \langle \Lambda^{(t_1,t_2)}, R^{(t_1,t_2)}  \rangle, \quad  \text{ for } (t_1,t_2) \in \tilde \Eomt  \label{eq:bca1}\\
		\lambda_t &\leftarrow \argmaxwrt[\lambda_t\in\RR^{n_t}]   - \epsilon \langle \bK,  \bU \rangle - \langle \lambda_t, \mu_t  \rangle,  \quad  \text{ for } t \in \tilde \Vomt_= \label{eq:bca2} \\
		\lambda_t &\leftarrow \argmaxwrt[\lambda_t\in\RR^{n_t}_+]   - \epsilon \langle \bK,  \bU \rangle - \langle \lambda_t, d_t  \rangle,  \quad  \text{ for } t \in \tilde \Vomt_\leq. \label{eq:bca3}
		\end{align}
	\end{subequations}
	The objectives of the unconstrained problems \eqref{eq:bca1} and \eqref{eq:bca2} are strictly concave, and thus a necessary and sufficient condition for optimality is that the respective gradient vanishes.
	Note that for each $(t_1,t_2)\in \Eomt$ the gradient of \eqref{eq:bca1} with respect to $\Lambda^{(t_1,t_2)}$ is
	\begin{equation}
	\exp(-\Lambda^{(t_1,t_2)}/\epsilon) \odot \left( \sum_{i_0,\dots,i_{\ccT} \setminus \{ i_{t_1}, i_{t_2} \} } \bK_{i_0\dots i_\ccT} \left(\prod_{t\in \tilde \Vomt}(u_t)_{i_t}\right) \left(\prod_{(\tau_1, \tau_2)\in \tilde \Eomt \setminus (t_1,t_2) } U_{i_{\tau_1}  i_{\tau_2}}^{(\tau_1, \tau_2)}\right) \right) - R^{(t_1,t_2)},
	\end{equation}
	and setting it to zero gives \eqref{eq:sinkhorn_graph_E}. Similarly, for $t\in\Vomt_=$ the gradient of \eqref{eq:bca2} with respect to $\lambda_t$ is 
	\begin{equation}
	\exp(-\lambda_t/\epsilon) \odot \left( \sum_{i_0,\dots,i_{t-1},i_{t+1},\dots,i_{\ccT} } \bK_{i_0\dots i_\ccT} \left(\prod_{\tau \in \tilde \Vomt \setminus \{t\} }(u_\tau)_{i_\tau}\right)  \left(\prod_{(t_1,t_2)\in \tilde \Eomt}U_{i_{t_1}  i_{t_2}}^{(t_1, t_2)}\right) \right) - \mu_t,
	\end{equation}
	which yields \eqref{eq:sinkhorn_graph_Veq}.
	Finally, note that the objective in \eqref{eq:bca3} can be written as
	\begin{equation}
	\sum_{i_t}  \left( -  \epsilon e^{-(\lambda_t)_{i_t} /\epsilon} \left( \sum_{\substack{i_0,\dots,i_{t-1}, i_{t+1},\dots,i_\ccT}} \bK_{i_0\dots,i_\ccT} \left(\prod_{\tau \in \tilde \Vomt \setminus \{t\} }(u_\tau)_{i_\tau}\right)  \left(\prod_{(t_1,t_2)\in \tilde \Eomt}U_{i_{t_1}  i_{t_2}}^{(t_1, t_2)}\right) \right)  -  (\lambda_t)_{i_t} d_{i_t} \right).
	\end{equation}
	Thus, the maximization in \eqref{eq:bca3} can be performed in each element of $\lambda_t$ individually. 
	If the derivative of the objective in \eqref{eq:bca3} with respect to $(\lambda_t)_{i_t}$ vanishes for a feasible, i.e., non-negative, point, then this is the global maximizer. Otherwise, the maximizer is the projection on the feasible set, i.e., $(\lambda_t)_{i_t}=0$.
	This yields \eqref{eq:sinkhorn_graph_Vineq}. The linear convergence of the scheme follows from \cite{luo1992convergence}.
	
In case only the assumptions in Theorem~\ref{thm:solution_extension} are satisfied, we perform a block coordinate ascent in the dual of \eqref{eq:omt_multi_graph_hat}. The dual variables can then be augmented by zero-entries as in \eqref{eq:u_augment} to arrive at the scheme in the Proposition. \hfill \Halmos
\endproof

Computing the projections of $\bK \odot \bU$ in Proposition~\ref{prp:sinkhorn} is in general still expensive, since computing the sums in \eqref{eq:proj_discrete} and \eqref{eq:proj_discrete_bi} requires $\mathcal{O}(n^\ccT)$ operations.
However, in the dynamic minimum-cost flow problems, there are additional structures in the cost tensor $\bC$, and thus in the tensor $\bK$. Namely, these tensors decouple according to the graphs in Figures~\ref{fig:path_graph} and \ref{fig:multicommodity}. The next subsections describe how these structures can be utilized in order to efficiently compute the projections needed to apply the scheme in Proposition~\ref{prp:sinkhorn}.

\subsection{Sinkhorn's method for the dynamic minimum-cost flow problem.}
\label{subsec:scheme_singlecommodity}

Recall that the dynamic minimum-cost flow problem \eqref{eq:singlecommodity} is a multi-marginal optimal transport problem. In particular, it can be written on the form \eqref{eq:omt_multi_graph}, where $ \tilde \Vomt_= = \{1,\ccT\}$, $\tilde \Vomt_\leq = \{2,\dots,\ccT-1\}$ and $\Eomt = \emptyset$.
Adding the entropy term \eqref{eq:entropy_term} yields then an entropy regularized problem \eqref{eq:omt_multi_graph_reg}, which in this case explicitly reads
\begin{equation} \label{eq:singlecommodity_reg} 
	\begin{aligned} 
	\min_{\bM\in \RR_+^{ n^\ccT}} \quad & \langle \bC, \bM\rangle  + \epsilon D(\bM) \\ 
	\mbox{subject to } \qquad 
	&P_{t}(\bM)\le d_t, \qquad \mbox{for } \qquad t=2,\ldots, \ccT-1 \\
	& P_{1}(\bM)=\mu_1, \\
	& P_{\ccT}(\bM)= \mu_t,
\end{aligned}
\end{equation}
where $\bC$ is defined by
\begin{equation*} %\label{eq:C_tensor_single_commodity}
\bC_{i_1\dots i_\ccT} = \sum_{t=1}^\ccT c_{i_t} + \sum_{t=1}^{\ccT-1} C_{i_t i_{t+1}}.
\end{equation*}

\begin{remark}
Without the inequality constraints $P_{t}(\bM)\le d_t$, and with zero cost on the edges, $c=\mathbf{0}$, the entropy-regularized problem \eqref{eq:singlecommodity_reg} is a discrete Schr\"odinger bridge problem \cite{pavon2010discrete, haasler19ensemble,haasler20trees}.
The Schr\"odinger brige problem is tightly connected to optimal transport \cite{chen2016relation, leonard2013schrodinger}.
It is a popular tool in ensemble control applications, as it provides a framework for steering a given distribution, i.e., an ensemble of agents, to a target one \cite{chen2016optimalPartI,brockett2012notes}.
In particular, network flow problems of this form have previously been considered in \cite{chen16network,CheGeoPavTan17,CheGeoPavTan19}.
This connection to the Schr\"odinger bridge problem gives another motivation for adding the regularizing entropy term to the objective of \eqref{eq:singlecommodity}. Namely, the Schr\"odinger bridge problem on a network can be interpreted as an ensemble of agents, which are each evolving according to a Markov chain \cite{haasler19ensemble, pavon2010discrete}. The entropy term thus induces a stochastic component to the problem, which yields a more smoothed out solution. Therefore, the solutions to the regularized problem \eqref{eq:singlecommodity_reg} can be understood as robust transport plans \cite{chen16network,CheGeoPavTan17,CheGeoPavTan19}.
\end{remark}

According to Theorem~\ref{thm:solution_extension}, the solution to the regularized problem \eqref{eq:singlecommodity_reg} is of the form $\bM= \bK \odot \bU$, where
\begin{equation}
\bK_{i_1 \dots i_\ccT} = \left( \prod_{t=1}^\ccT k_{i_t} \right) \left( \prod_{t=1}^{\ccT-1} K_{i_t i_{t+1}} \right),
\end{equation}
with $k=\exp(-c/\epsilon)$ and $K=\exp(-C/\epsilon)$, and $\bU= u_1 \otimes \dots \otimes u_\ccT$. The components of the tensor $\bU$ can be found utilizing Proposition~\ref{prp:sinkhorn}. 
In particular, the solution is found by iterating
\begin{equation} \label{eq:sinkhorn_singlecommodity}
\begin{aligned}
u_t &\leftarrow u_t \odot \mu_t ./ P_{t} (\bK \odot \bU), \qquad \qquad \quad \; \, \text{ for } t=1,\ccT, \\
u_t &\leftarrow \min \left(  u_t \odot d ./ \left( P_t(\bK \odot \bU)\right), \ett  \right), \quad \text{ for }   t =2,\dots, \ccT-1. 
\end{aligned}
\end{equation}
In this case, where the cost decouples according to a path graph, the projections can be computed efficiently \cite[Proposition 2]{elvander19multi}.
Namely, the projections for this problem are of the form
\begin{equation} \label{eq:proj_bridge}
P_t(\bK \odot \bU) = u_t \odot k_t \odot \hat \varphi_t \odot \varphi_t,
\end{equation}
for $t=1,\dots,\ccT$, where
\begin{subequations}\label{eq:phi_phi_hat}
\begin{align}
\hat \varphi_t & = K^T \diag(u_{t-1} \odot k_{t-1}) K^T \dots \diag(u_2 \odot k_2) K^T (u_1 \odot k_1), \label{eq:phi_hat} \\
\varphi_t & = K \diag(u_{t+1} \odot k_{t+1}) K \dots \diag( u_{\ccT-1} \odot k_{\ccT-1}) K (u_\ccT \odot k_\ccT). \label{eq:phi}
\end{align}
\end{subequations}
The Sinkhorn algorithm \eqref{eq:sinkhorn_singlecommodity}  is summarized in Algorithm~\ref{alg:single_commodity}.
\begin{algorithm}[tb]
	\begin{algorithmic}
		\STATE Initialize $u_1,\dots,u_\ccT$, $t=1$, $\hat \varphi_1=\ett$, $\varphi_\ccT=\ett$
		\WHILE{Not converged}	
		\FOR{ $t=\ccT-1,\dots,1$}
		\STATE Update $ \varphi_t \leftarrow K (u_{t+1} \odot k_{t+1} \odot \varphi_{t+1})$ 
		\ENDFOR
		\STATE Update $u_1 \leftarrow \mu_1./ \varphi_1$
		\FOR{$t=2,\dots,\ccT-1$}
		\STATE Update $ \hat \varphi_t \leftarrow K^T (u_{t-1} \odot k_{t-1} \odot \hat \varphi_{t-1})$ 
		\STATE Update $u_t \leftarrow \min \left( d_t ./ (\varphi_t \odot \hat \varphi_t \odot k_t) , \ett \right)$
		\ENDFOR
		\STATE Update $ \hat \varphi_\ccT \leftarrow  K^T (u_{\ccT-1} \odot k_{\ccT-1} \odot \hat \varphi_{\ccT-1}) $
		\STATE Update $u_\ccT \leftarrow \mu_\ccT ./ \hat\varphi_\ccT$
		\ENDWHILE
		\RETURN $u_1,\dots,u_\ccT$
	\end{algorithmic}
	\caption{ Scheme for solving the dual of the regularized dynamic flow problem \eqref{eq:singlecommodity_reg}. }
	\label{alg:single_commodity}
\end{algorithm}

Note that intermediate results of \eqref{eq:phi_hat} and \eqref{eq:phi} are stored, and the updates in \eqref{eq:sinkhorn_singlecommodity} are scheduled such that for each update only one matrix-vector multiplication needs to be performed. Thus, in the case of a dense matrix $K$, one iteration sweep, i.e., once updating all vectors $u_t$, for $t=1,\dots,\ccT$, is of complexity $\ccO(\ccT n^2)$.
However, for sparse networks the matrix $K$ is also sparse, and thus the matrix multiplications required to compute the projections \eqref{eq:proj_bridge} via \eqref{eq:phi_phi_hat} become even more efficient, as discussed in the following remark.
\begin{remark}\label{rem:computational_complexity}
	Note that $K_{ij}=0$ if $(i,j) \notin \Eflow$, and $K_{ij}>0$ if $(i,j)\in \Eflow$. Thus, multiplication with a vector $v\in \RR^n$ can be performed as
	\begin{equation}
	\left( Kv \right)_i = \sum_{j \in N(i)} K_{ij} v_j.
	\end{equation}
	This multiplication is of order $\mathcal{O}(\Delta(\Nflow) \cdot n)$, where $\Delta(\Nflow)$ is the maximum degree of $\Nflow$, i.e., the highest number of neighboring nodes among the nodes $\Vflow$.
	The complexity of one iteration sweep in Algorithm~\ref{alg:single_commodity} is thus $\ccO(\ccT n \Delta(\Nflow))$.
\end{remark}

\subsection{Sinkhorn's method for the dynamic multi-commodity minimum-cost flow problem.}
\label{subsec:scheme_multicommodity}

Similarly to the previous section, the multi-commodity problem \eqref{eq:multicommodity} is also a multi-marginal optimal transport problem of the form \eqref{eq:omt_multi_graph}.
In particular, here the constraint sets are $\tilde \Vomt_= = \emptyset$, $\tilde \Vomt_\leq = \{ 2,\dots, \ccT-1\}$ and $\tilde \Eomt=\{(0,1), (0,\ccT)\}$.
Regularizing the problem with an entropy term, it is of the form \eqref{eq:omt_multi_graph_reg}, which in this case reads
\begin{equation} \label{eq:multicommodity_reg} 
\begin{aligned}
\minwrt[\bM\in \RR_+^{L\times n^\ccT}] \ & \  \langle \bC, \bM \rangle  + \epsilon D(\bM)  \\ 
\mbox{subject to } \
&\ P_{0,1}(\bM)=R^{(0,1)},\\
& \ P_{0,\ccT}(\bM)=R^{(0,\ccT)}, \\
&\ P_{t}(\bM)\le d_t, \qquad \qquad \; \mbox{ for } t=2,\ldots, \ccT-1, 
\end{aligned}
\end{equation}
where $\bC$ is defined by
\begin{equation}
\bC_{i_0\dots i_\ccT} = \sum_{t=2}^{\ccT-1} (\CL)_{i_0 i_t} +  \sum_{t=1}^{\ccT-1} C_{i_t i_{t+1}}.
\end{equation}
The solution to \eqref{eq:multicommodity_reg} can again be expressed in terms of its dual variables, as described in Theorem~\ref{thm:solution_extension}. In particular, the optimal mass transport plan is of the form $\bM=\bK \odot \bU$, where $\bK$ factorizes as
\begin{equation}\label{eq:K_multicommodity}
\bK_{i_0\dots i_\ccT} = \left( \prod_{t=2}^{\ccT-1} (\KL)_{i_0 i_t} \right) \left(  \prod_{t=1}^{\ccT-1} K_{i_t i_{t+1}} \right),
\end{equation}
where $\KL=\exp(-\CL/\epsilon)$ and $K=\exp(-C/\epsilon)$.
Moreover, the tensor $\bU$ is of the form
\begin{equation} \label{eq:U_multicommodity}
\bU_{i_0 \dots i_\ccT} = U^{(0,1)}_{i_0 i_1} U^{(0,\ccT)}_{i_0 i_\ccT} \prod_{t=2}^{\ccT-1} (u_t)_{i_t},
\end{equation}
and its components can be found according to Propositions~\ref{prp:sinkhorn} by iteratively updating
\begin{equation} \label{eq:sinkhorn_multicommodity}
\begin{aligned}
U^{(0,t)} & \leftarrow U^{(0,t)} \odot R^{(0,t)} ./ P_{0,t} (\bK \odot \bU), \quad \text{ for } t=1,\ccT, \\
u_t & \leftarrow \min \left( u_t \odot d ./ P_t(\bK \odot \bU), \ett \right), \quad \text{ for } t=2,\dots,\ccT-1.
\end{aligned}		
\end{equation}

Again, the tensor $\bK \odot \bU$ has a graph structure, which is illustrated in Figure~\ref{fig:multicommodity}. This graph contains cycles, and thus the results from \cite{haasler20trees} cannot be utilized.
Nevertheless, the projections can be computed relatively efficiently, as demonstrated by the next theorem.
\begin{theorem} \label{thm:biproj}
	Consider the tensors $\bK=\exp(-\bC/\epsilon)$, with $\bC$ defined as in \eqref{eq:cost_multicommodity} and $\epsilon>0$, and $\bU$ in \eqref{eq:U_multicommodity}. 
	With the matrices $\KL=\exp(-\CL/\epsilon)$ and $K=\exp(-C/\epsilon)$, define
	\begin{equation} \label{eq:psi_hat}
	\hat \Psi_t = \begin{cases} %U^{0,1} , \qquad t=1,\\
	 U^{(0,1)}  K , \qquad t=2,\\
	\left( \hat \Psi_{t-1} \odot \KL \right) \diag(u_{t-1}) K  , \qquad t= 3,\dots,\ccT,
	\end{cases} 
	\end{equation}
	and
	\begin{equation} \label{eq:psi}
	\Psi_t = \begin{cases} % U^{(0,\ccT)} , \qquad t=\ccT,\\
	U^{0,\ccT}  K^T , \qquad t=\ccT-1,\\
	\left( \Psi_{t+1} \odot \KL \right) \diag(u_{t+1}) K^T , \qquad t=1,\dots, \ccT-2.
	\end{cases}
	\end{equation}
	Then, the bi-marginal projections of the tensor $\bK \odot \bU$ are
	\begin{equation} \label{eq:proj_0t}
	\begin{aligned}
	P_{0,1} (\bK \odot \bU) & = U^{(0,1)} \odot  \Psi_1 \\
	P_{0,\ccT} (\bK \odot \bU) & = U^{(0,\ccT)} \odot \hat \Psi_\ccT \\
		P_{0,t} (\bK \odot \bU) & = \left(   \hat \Psi_{t} \odot   \Psi_t \odot \KL \right) \diag(u_t) , \quad \text{ for } t=2,\dots,\ccT-1.  
	\end{aligned}
	\end{equation}
\end{theorem}

\proof{Proof:}
Note that the tensor $\bK=\exp(-\bC/\epsilon)$ is element-wise defined as in \eqref{eq:K_multicommodity}, thus the bi-marginal projections of the tensor $\bK \odot \bU$ on the marginals $0$ and $t$, where $t \in \{2,\dots,\ccT-1\}$, are given by
\begin{equation}
\begin{aligned}
P_{0,t} (\bK \odot \bU) =& \sum_{\substack{i_1,\dots, i_{t-1} \\ i_{t+1},\dots, i_\ccT}} \left( \prod_{s=1}^{\ccT-1} K_{i_s i_{s+1}} \right) \left( \prod_{s=2}^{\ccT-1} (\KL)_{i_0 i_s} \right) U^{(0,1)}_{i_0 i_1} U^{(0,\ccT)}_{i_0 i_\ccT} \prod_{s=2}^{\ccT-2} (u_s)_{i_s} \\
=&  (u_t)_{i_t} (\KL)_{i_0 i_t} (\hat \Psi_{t})_{i_0 i_t} (\Psi_{t})_{i_0 i_t},
\end{aligned}
\end{equation}
where
\begin{equation}
\hat \Psi_t = \sum_{i_1,\dots, i_{t-1}} U^{(0,1)}_{i_0 i_1} K_{i_1 i_2} \left( \prod_{s=2}^{t-1} ( \KL \diag(u_s) )_{i_0 i_s} K_{i_s i_{s+1}} \right),
\end{equation}
and
\begin{equation}
\Psi_t = \sum_{i_{t+1},\dots, i_{\ccT}} U^{(0,\ccT)}_{i_0 i_\ccT} K_{i_{\ccT-1} i_\ccT} \left( \prod_{s=t+1}^{\ccT-1} ( \KL \diag(u_s) )_{i_0 i_s} K_{i_{s-1} i_{s}} \right).
\end{equation}
These terms lead to the recursive definitions of $\hat \Psi_t$ and $\Psi_{t}$ in \eqref{eq:psi_hat} and \eqref{eq:psi}. The projections $P_{0,1}(\bK \odot \bU)$ and $P_{0,\ccT}(\bK \odot \bU)$ are derived similarly. \hfill \Halmos
\endproof

The projections on one marginal can then be found by projecting the bi-marginal projections in \eqref{eq:proj_0t} on one of the marginals, which yields the following.
\begin{corollary} \label{cor:proj}
	The marginals of the tensor $\bK \odot \bU$ in Theorem~\ref{thm:biproj} are given by
	\begin{equation}
	\begin{aligned}
	P_t(\bK \odot \bU) &= u_t \odot  \left(   \hat \Psi_{t} \odot   \Psi_t \odot \KL  \right)^T \ett, \quad \text{for } t=2,\dots, \ccT-1, \\
	P_0(\bK \odot \bU) &=  \left(   \hat \Psi_{t} \odot   \Psi_t \odot \KL  \right) u_{t}.
	\end{aligned}
	\end{equation}
\end{corollary}

Theorem~\ref{thm:biproj} and Corollary~\ref{cor:proj} describe an efficient way to compute the projections required for the Sinkhorn scheme \eqref{eq:sinkhorn_multicommodity}, and the resulting computational method is summarized in Algorithm~\ref{alg:multi_commodity}.
\begin{algorithm}[tb]
	\begin{algorithmic}
		\STATE Initialize $u_2,\dots,u_{\ccT-1}$, $U^{(0,1)}$, $U^{(0,\ccT)}$
		\STATE Compute $\Psi_t$, for $t=1,\dots,\ccT$
		\WHILE{Not converged}	
		\STATE Update $U^{(0,1)} \leftarrow R^{(0,1)}./ \Psi_1$
 		\STATE Update $\hat \Psi_2 \leftarrow U^{(0,1)} K$ 
		\FOR{ $t=2,\dots,\ccT-1$}
		\STATE Update $u_t \leftarrow \min\left( d ./ ( (   \hat \Psi_{t} \odot   \Psi_t \odot K )^T \ett) \ ,\  \ett \right)$
		\STATE Update $\hat \Psi_{t+1} \leftarrow ( \hat \Psi_{t} \odot \KL ) \diag(u_{t}) K $ 
		\ENDFOR 
		\STATE $U^{(0,\ccT)} \leftarrow R^{(0,\ccT)}./ \hat \Psi_\ccT$
		\STATE Update $ \Psi_{\ccT-1} \leftarrow U^{(0,\ccT)} K^T$ 
		\FOR{ $t=\ccT-1,\dots, 2$}
		\STATE  Update $ \Psi_{t-1} \leftarrow ( \Psi_{t} \odot \KL ) \diag(u_{t}) K^T $
		\ENDFOR 
		\ENDWHILE
		\RETURN $u_2,\dots,u_{\ccT-1}$, $U^{(0,1)}$, $U^{(0,\ccT)}$
	\end{algorithmic}
	\caption{Scheme for solving the dual of the regularized dynamic multi-commodity flow problem \eqref{eq:multicommodity_reg}
	}
	\label{alg:multi_commodity}
\end{algorithm}
Similarly to the algorithm for the single-commodity setting, intermediate results can be stored and utilized.

\begin{remark} 
	The computational bottleneck of the Sinkhorn iterations lies in computing the projections.
One iteration sweep of the Sinkhorn iterations requires updating each of the matrices in \eqref{eq:psi_hat} and \eqref{eq:psi} once.
For dense matrices $K$ each of these updates is of complexity $\ccO(Ln^2)$, and thus one full iteration sweep can be done in $\ccO( \ccT Ln^2)$.
	However, as noted in Remark~\ref{rem:computational_complexity}, the matrix $K$ inherits the sparsity of the network, and this can be exploited to perform the matrix multiplications in \eqref{eq:psi_hat} and \eqref{eq:psi} more efficiently.
Thus, the complexity of the matrix-matrix multiplication is decreased to $\mathcal{O}(\Delta(\Nflow) \cdot L n)$, and one full iteration sweep can be done in $\ccO( \ccT \Delta(\Nflow) Ln)$.

\end{remark}

\begin{remark} \label{rem:multi_tensor_scheme}
	In Section~\ref{sec:mc_form} we have formulated the multi-tensor problem \eqref{eq:multicommodity_tensor_problem} as the one-tensor problem \eqref{eq:multicommodity} in order to bring it on the form of a graph-structured optimal transport problem \eqref{eq:omt_multi_graph} and then solve it.
	Alternatively, we could have regularized each of the $L$ optimal transport problems in \eqref{eq:multicommodity_tensor_problem} separately, yielding the regularized problem
\begin{equation} \label{eq:multicommodity_tensor_problem_reg}
\begin{aligned} 
\minwrt[ \bM^1,\dots,\bM^L \in \RR_{+}^{n^{\ccT}}  ] \ &  \sum_{\ell=1}^L \left( \langle \bC, \bM^\ell \rangle + \epsilon D(\bM^\ell) \right) \\
\text{subject to } \ & \ P_{1} ( \bM^\ell ) = \mu_1^\ell, \quad \ell=1,\dots,L, \\
& \ P_{\ccT} ( \bM^\ell ) = \mu_\ccT^\ell, \quad \ell=1,\dots,L, \\
& \ \sum_{\ell=1}^L P_t( \bM^\ell) \leq d \quad t=2,\dots,\ccT-1,
\end{aligned}
\end{equation}
where $\bC$ is defined as in \eqref{eq:T_cost}.
In fact, this problem is equivalent to the regularized problem \eqref{eq:multicommodity_reg}.
Moreover, in this representation the Sinkhorn iterations are given by
\begin{equation} \label{eq:Sinkhorn_multi_tensor}
\begin{aligned}
u_1^\ell &\leftarrow u_1^\ell \odot \mu_1^\ell ./ P_{1} (\bK^\ell \odot \bU^\ell), \quad &&\text{ for }  \ell =1,\dots, L \\
u_t &\leftarrow \min \left(  u_t \odot d ./ \left( \sum_{\ell\in L} P_t(\bK^\ell \odot \bU^\ell)\right), \ett  \right), \quad &&\text{ for }   t =2,\dots, \ccT-1 \\
u_\ccT^\ell &\leftarrow   u_\ccT^L \odot \mu_\ccT^\ell ./ (P_\ccT(\bK^\ell \odot \bU^\ell)), \quad &&\text{ for } \ell =1,\dots, L
\end{aligned}
\end{equation}
 and these are equivalent to the Sinkhorn iterations derived above (cf. \eqref{eq:sinkhorn_multicommodity}).
Recall from Section~\ref{subsec:discussion} that one convenient feature of formulation \eqref{eq:multicommodity_tensor_problem} is that it can be easily extended to allow for commodities that enter and leave the network at different times.
Therefore, as can be seen here, such problems can also be solved efficiently.
\end{remark}

\section{Simulations.}
\label{sec:siumlations}

In this Section we illustrate the computational efficiency of our propsed framework. First, we compare its performance with a Simplex solver on two different types of networks.
Finally, we illustrate it in a traffic routing problem with a large number of commodities.

\subsection{Performance study on a sparse grid network.} \label{subsec:sparse}

We first consider a dynamic multi-commodity minimum-cost network flow problem on a sparse network.
To this end, let $\Nflow$ be a grid of $10\times 10$ nodes, and let the source $\ccS^+$ for all commodities be an incoming edge to one corner of the square, and let the sink $\ccS^-$ be an outgoing edge from the opposite corner. Thus, the total number of directed edges is $n=84$. Moreover, in this set-up the sink and source can be understood as the two corner vertices.
We consider the case of $L=50$ commodities, and let the total flow of each commodity be $1$, that is $\mu_0= \ett$. Moreover, the capacity vector $d\in \RR_{+}^n$ is defined as $d_i=L$ for $i \in \{ \ccS^+ \cup \ccS^- \}$, and $d_i=1$ otherwise. 
Here we do not allow for intermediate storage on the vertices or the edges, except in the sink and source.
This problem is solved for a time horizon of $\ccT=80$ utilizing Algorithm~\ref{alg:multi_commodity}.
We also solve the problem in node-edge formulation (cf. \cite{ford1962flows, tomlin1966minimum}) in the time-expanded network using the solver CPLEX \cite{cplex2009v12}. In our experiments we observed that CPLEX performs best when using the dual simplex algorithm, and we thus assign this algorithm when calling CPLEX to decrease its start-up time.
We run this problem for $10$ different experiments, where in each trial the cost for a unit flow of each commodity on each edge is randomly assigned from a uniform distribution on $[0,1]$, that is we let $c_{e}^\ell \sim \text{Unif}([0,1])$, for $\ell=1,\dots,L$, and $e\in \Eflow$.

\begin{figure}
\centering
\includegraphics[width=\textwidth]{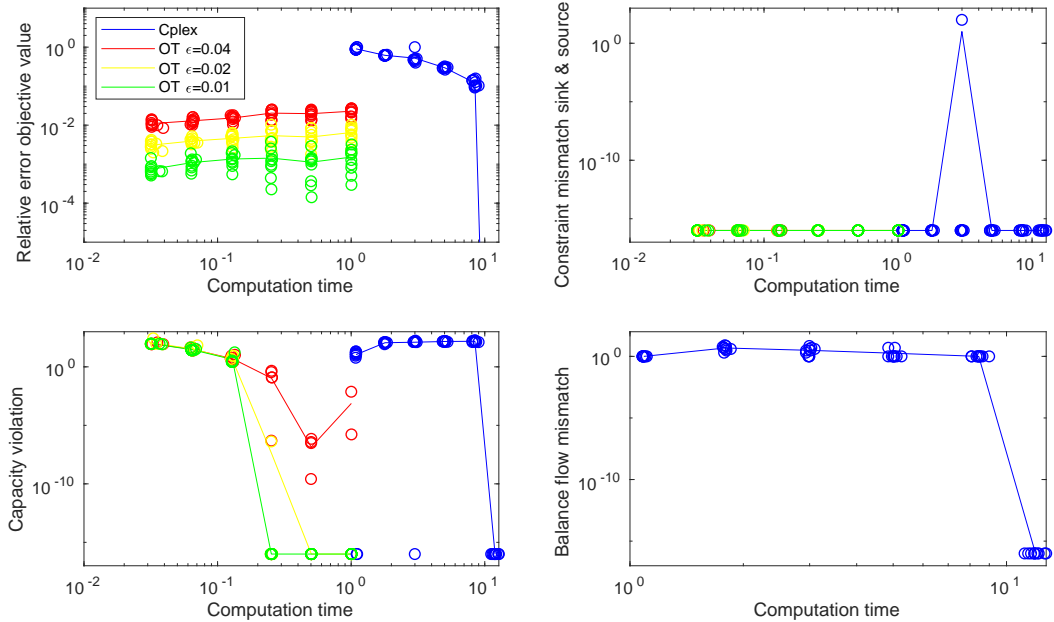}
\caption{Performance study on a sparse grid network. The plots show the $L_1$-norm of the objective value error and constraint mismatches of the current solution over time. Circles denote each experiments outcome, and the full lines are plotted between their means. The machine precision is $10^{-16}$. }
\label{fig:performance_sparse}
\end{figure}

Figure~\ref{fig:performance_sparse} shows some measures of error as a function of computation time.
The error in the objective value compares the objective value of the current solution to the optimal objective value. Note that our proposed algorithm is based on the regularized problem \eqref{eq:multicommodity_reg}, and therefore cannot achieve the true objective value. However, the smaller the regularization parameter $\epsilon$, the closer we get to the true optimum.
Since CPLEX utilizes a dual simplex method, its solution becomes meaningful only after the last iteration, in the sense that not all constraints are fulfilled for the intermediate iterates. In particular, this is the case for the flow balance constraint. In contrast, when solving \eqref{eq:multicommodity_reg} using Algorithm~\ref{alg:multi_commodity}, the intermediate iterates by construction satisfy the flow balance constraints in the nodes, and they also satisfy the mismatch in the sinks and sources to machine precision. Moreover, Algorithm~\ref{alg:multi_commodity} converges linearly to the optimal solution of \eqref{eq:multicommodity_reg}.
With the smallest tested regularization parameter, $\epsilon=0.01$, the capacity constraints are satisfied in about $0.1$ seconds.
On the other hand CPLEX takes more than $10$ seconds to find a solution, which is two orders of magnitude longer than Algorithm~\ref{alg:multi_commodity}.
Note that state-of-the-art methods for multi-commodity flows can typically not be expected to improve the run time by more than an order of magnitude as compared to standard LP solvers \cite{barnhart2009multicommodity, retvdri2004novel, khodayifar2019minimum}.
Our proposed algorithm is thus competitive with specialized state-of-the-art methods for network flow problems.

\subsection{Performance study on a dense random network.}

Next, we study the performance of Algorithm~\ref{alg:multi_commodity} in a less favourable setting. Here we consider a dense random network with $40$ nodes.
Between each (ordered) pair of nodes we create a directed edge with probability $1/2$.
The expected value of the number of edges in the network is thus $\binom{40}{2} = 780$.
Moreover, we allow for intermediate storage in the nodes, but not in the edges.
Therefore, we augment the state space by the set of nodes as described in Section~\ref{subsec:discussion}, and the expected size of the distributions support is thus $\mathbb{E}[n] = 820$.
We equip each of the $L=100$ commodities with a random source and sink on the set of nodes. The total flow of each commodity is set to $1$, i.e., $\mu_0 = \ett$, and the capacity vector $d\in\RR_{+}^n$ is defined as $d_i=1$, if $i\in \Eflow$, and $d_i=L$, if $i\in \Vflow$.
As in the previous example, the cost for each commodity and each edge is assigned from a uniform distribution on $[0,1]$. Moreover, the cost for intermediate storage on the nodes is $0$.
We solve the problem for $\ccT=100$ time intervals using Algorithm~\ref{alg:multi_commodity} and solve its node-edge formulation in the time-expanded network with the dual simplex algorithm in CPLEX. Here, the time expanded network has $4000$ nodes and in the mean $81180$ edges.

\begin{figure}
	\centering
	\includegraphics[width=\textwidth]{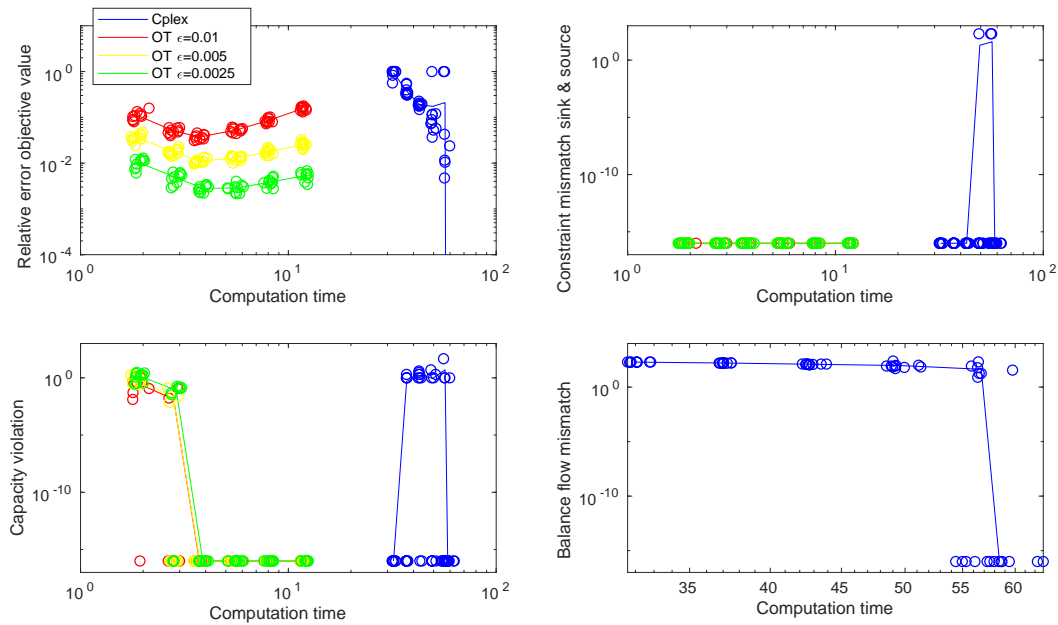}
	\caption{Performance study on a larger dense random network. The plots show the $L_1$-norm of the objective value error and constraint mismatches of the current solution over time. Circles denote each experiments outcome, and the full lines are plotted between their means. The machine precision is $10^{-16}$. } \label{fig:performance_dense}
\end{figure}

The performance for $10$ trials of the described setup is illustrated in Figure~\ref{fig:performance_dense}.
Qualitatively, we see a similar behavior as in the performance study for the sparse network in Section~\ref{subsec:sparse}.
In particular, as in the previous example, in contrast to the intermediate iterates produced by our method, the intermediate iterates produced by CPLEX do not correspond to flows since the flow balance constraint is in general not fulfilled.
Moreover, our method converges linearly to an optimal solution of \eqref{eq:multicommodity_reg}, and with the smallest tested regularization parameter $\epsilon=0.0025$ Algorithm~\ref{alg:multi_commodity} gives a good approximation to the optimal solution in about one second, whereas the CPLEX solver requires about 15 seconds to find a solution.
Even in the less favourable setting of a dense network with intermediate storage on the nodes we thus still get an aproximate solution in less than an order of magnitude of CPLEX's run-time.

\subsection{Traffic routing problem with a large number of commodities.}

\begin{figure}[h!] 
	\centering
	\includegraphics[width=0.5\textwidth]{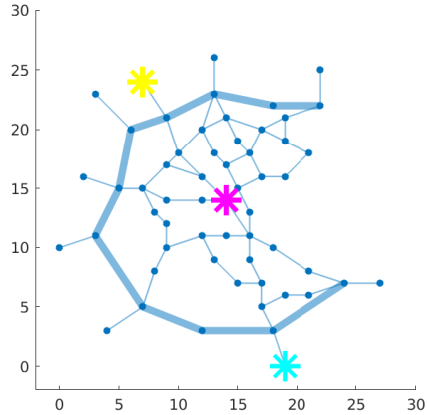}
	\caption{Map of a street network. Every edge represents two directed edges, one in each direction. Broader edges represent highways. The three stars represent three different commodities' sinks. }
	\label{fig:sthlm_map}
\end{figure}

\begin{figure} 
	\centering
	\includegraphics[width=\textwidth]{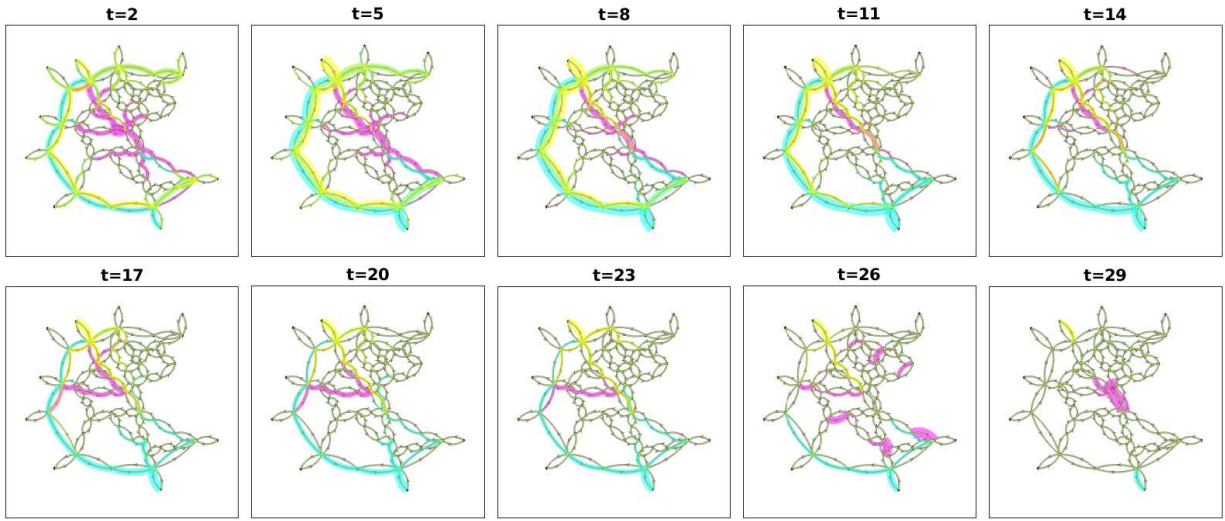}
	\caption{The optimal traffic flow over time for three of the commodities.}
	\label{fig:traffic_3comm}
\end{figure}

\begin{figure}
	\centering
	\subfigure[\ Number of agents in sources, roads, and sinks, over time for the three commodities in Figure~\ref{fig:traffic_3comm}. ]{\includegraphics[width=0.45\textwidth]{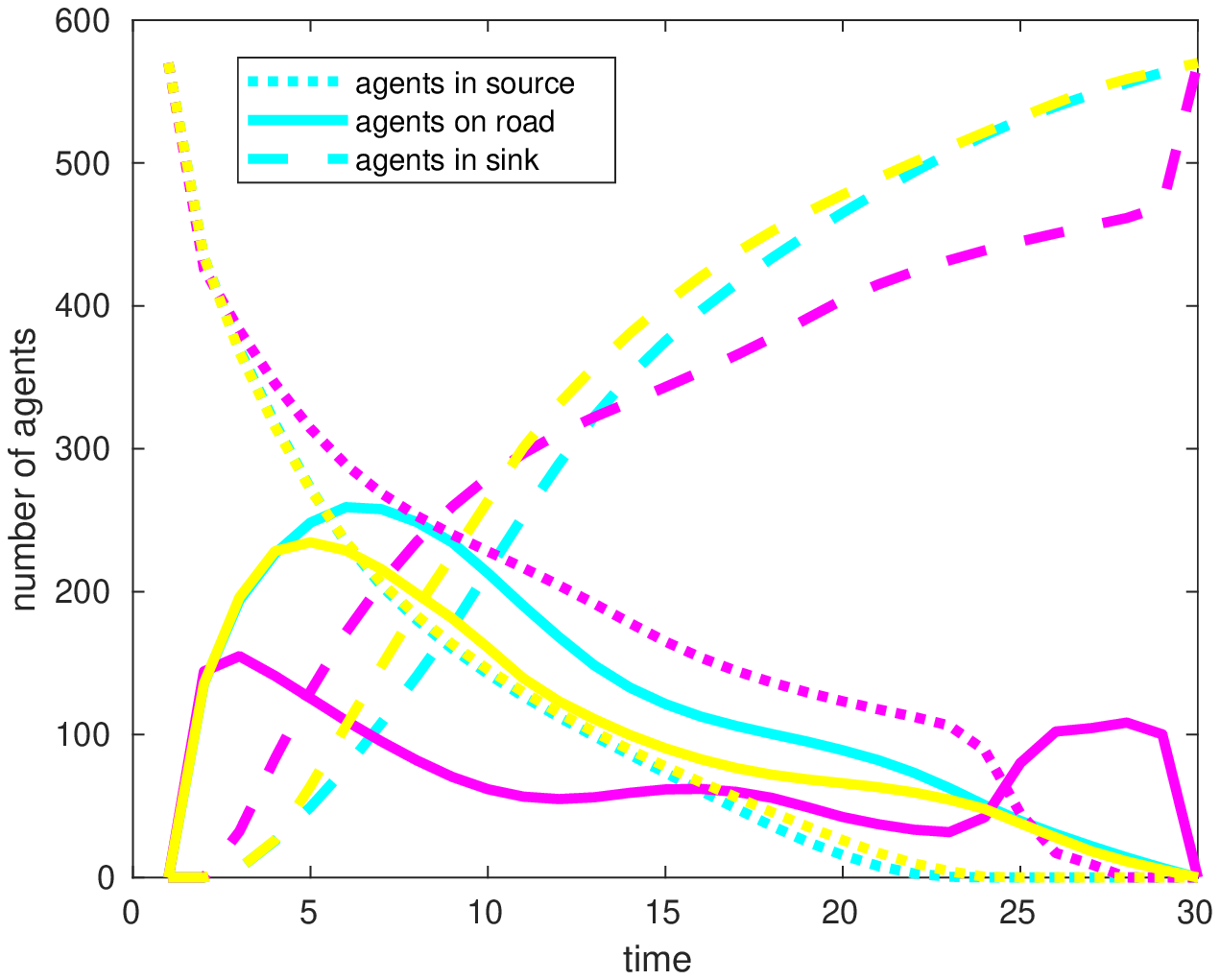} \label{fig:status_3comm} }
	\subfigure[\ Blue curves correspond to all agents in the scenario in Figure~\ref{fig:traffic_3comm}. Green and red curves describe the scenario, where the cost for staying in a source is equal to the cost on the roads (0.1) and to the cost for staying in a sink (0), respectively. ]{\includegraphics[width=0.45\textwidth]{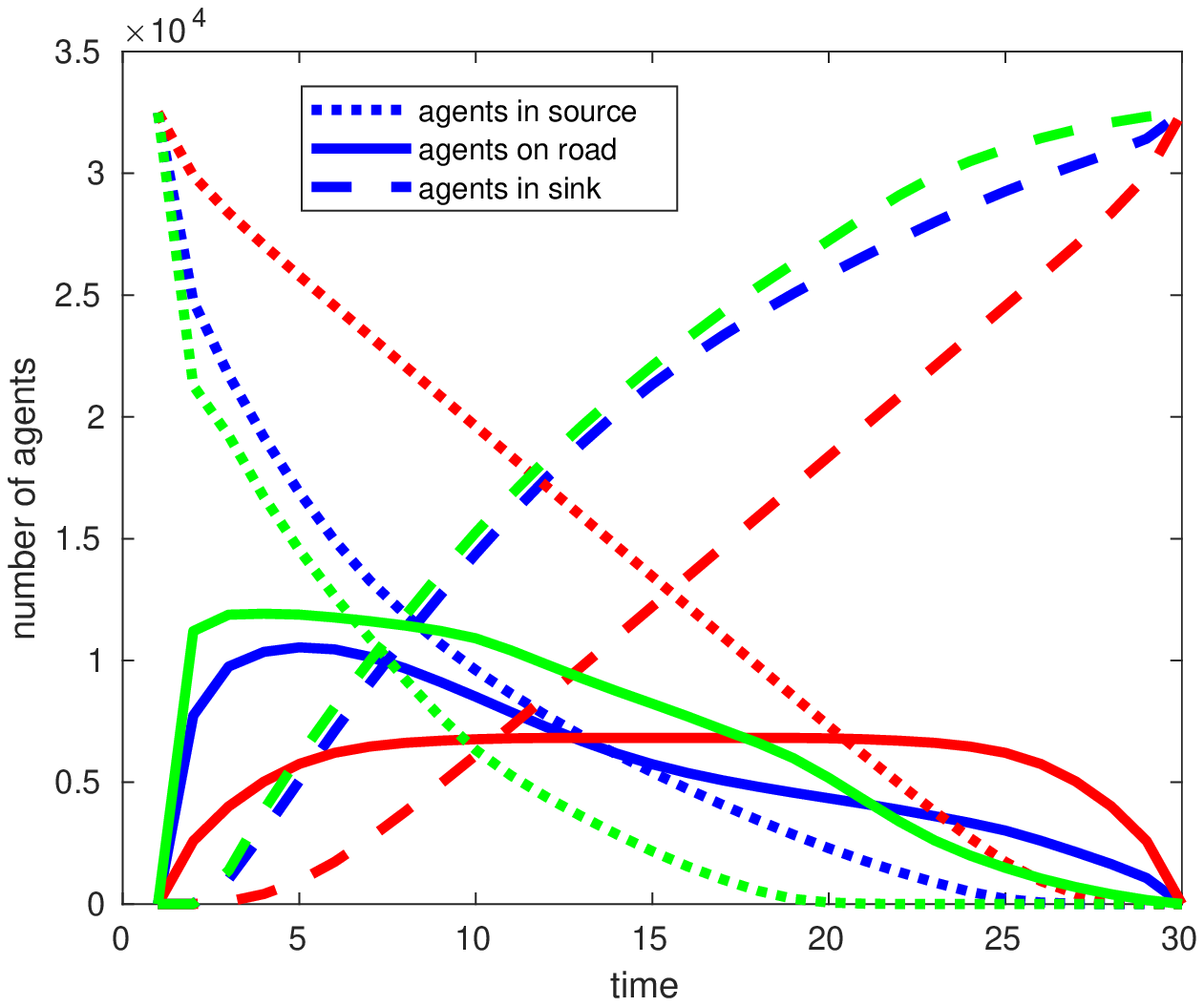} \label{fig:status_3comm_costs} }
	\caption{Agents status over time.}
\end{figure}

We apply our framework to a traffic routing problem in the street network illustrated in Figure~\ref{fig:sthlm_map}, which consists of $57$ nodes and $150$ directed edges.
Let every node be both a sink and a source, and as described in Section~\ref{subsec:discussion} we thus let the state space be of size $n= 150+ 2\cdot 57 = 264$.
Assume that there is an equal amount of $10$ agents travelling between every pair of nodes.
This can be modelled by associating each node with one commodity, and imposing that every commodity is initially uniformly distributed on the set of sources, and finally concentrated in the associated sink node. 
In particular, this means that the number of commodities is $L= 57$, and the two matrix constraints in \eqref{eq:multicommodity} are defined by the matrices $R^{(0,1)}, R^{(0,\ccT)} \in \RR_{+}^{L \times n}$ with entries
\begin {equation}
R^{(0,1)}_{\ell,i} = \begin{cases} 10 , & \text{ if }i \in \ccS_\ell^+ = \Vflow, \\
0 , & \text{otherwise,} \end{cases} \qquad
R^{(0,\ccT)}_{\ell,i} = \begin{cases} 570 , & \text{ if }i \in \ccS_\ell^-, \\
0  , & \text{otherwise.} \end{cases}
\end{equation}
We consider the scenario with intermediate storage in the edges, but without storage on the nodes. However, agents are permitted to stay in their respective sink and source, but once they leave their source they may not return to it, and once they reach their sink they may not leave it. This structure is imposed by the cost matrix $C$ as defined in \eqref{eq:cost_traffic_routing}. 
The wider streets in Figure~\ref{fig:sthlm_map} describe highways, and we denote the set of highways as $\ccH$.
Since our framework assumes uniform travel time on all edges, the fact that the roads in $\ccH$ are longer than the other roads models that agents can drive faster on the highway.
Let $l_i$ denote the Euclidean length of road $i\in\Eflow$.
We define the capacities for each state as
%edge as $d_i= 100 l_i$ if $i \in \ccH$, and $d_i = 20 l_i $ if $i\in \Eflow \setminus \ccH$. The capacities in the sinks and sources are $d_i = 100 L$, that is essentially unlimited.
\begin{equation}
d_i = \begin{cases} 100 l_i ,  & \text{ if } i \in \ccH , \\
20 l_i, &  \text{ if } i \in \Eflow \setminus \ccH , \\
 100L, & \text{ if } i \in \ccS . \end{cases}
\end{equation}
The cost for an agent to be in any of the states is defined in the matrix $C_L$. The costs are assumed equal for all agents and defined for all commodities $\ell=1,\dots, L$ as
\begin{equation}\label{eq:cost_trafic_example}
(C_L)_{\ell i} = \begin{cases} 0.01,  & \text{ if } i \in \ccS^+ , \\
 0.1, &  \text{ if } i \in \Eflow , \\
 0, & \text{ if } i \in \ccS^- . \end{cases}
\end{equation}
Thus, the central controller aims to minimize the time agents spend inside the network, and makes them reach the sink early rather than wait in the source.
We consider the problem with final time $\ccT=30$.
The problem is solved using Algorithm~\ref{alg:multi_commodity} with regularization parameter $\epsilon = 0.01$.
For the three commodities associated with the sinks highlighted in Figure~\ref{fig:sthlm_map}, the optimal flows are visualized in Figure~\ref{fig:traffic_3comm}.
One can see that traffic is sent at all places in the network, and finally concentrates towards the three sinks.
For the three commodities the amount of agents in the sources, roads, and sinks, respectively, is plotted over time in Figure~\ref{fig:status_3comm}.
The total flows distribution over source, roads, and sink over time can be seen in the blue lines in Figure~\ref{fig:status_3comm_costs}.
At the first time instance many agents are sent from the sources into the network.
Towards the end of the time interval less and less agents are on the roads.

We also vary the cost for agents to stay in the source, see Figure~\ref{fig:status_3comm_costs}.
Clearly, if the cost for being in a source is increased to $(C_L)_{\ell i} = 0.1$, for $i\in \ccS^+$ and $\ell=1,\dots,L$, more agents are sent into the network early on.
If the cost for being in a source is equal to being in a sink, i.e., $(C_L)_{\ell i} = 0$, for $i\in \ccS^+$ and $\ell=1,\dots,L$, the amount of flow on the roads over time looks very symmetric.

\begin{figure} 
	\centering
	\includegraphics[width=\textwidth]{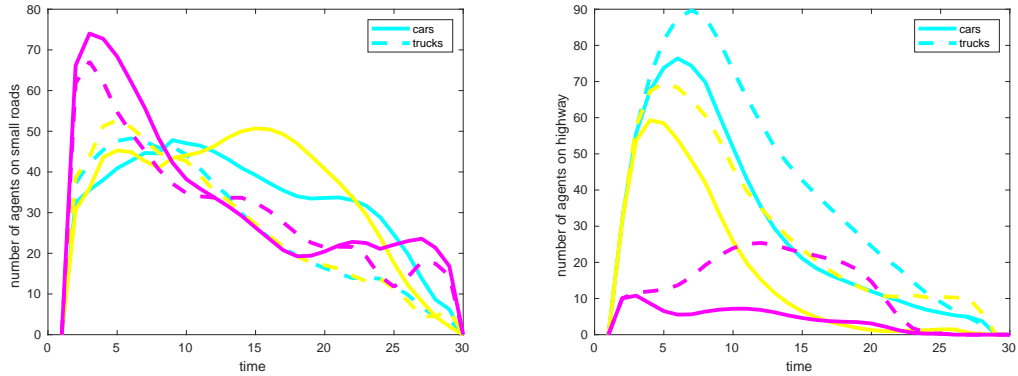}
	\caption{Distribution of six commodities on small roads and highways over time.}
	\label{fig:status_3comm_trucks_roads}
\end{figure}

Finally, we consider a scenario where a second type of commodity is present in the network.
Therefore, the total amount of commodities is increased to $2L=114$.
We interpret the first set of $L$ commodities as cars and denote them as $\ccL_C$.
The second set of $L$ commodities are interpreted as trucks and denoted by $\ccL_T$.
For each set of commodities, the initial and final distributions are defined as before, but the number of agents in each commodity is halved in order to get the same total number of agents .
That is, we define the new constraint matrices $\hat R^{(0,1)}, \hat R^{(0,\ccT)} \in \RR_{+}^{2L \times n}$ as 
\begin{equation}
\hat R^{(0,1)} = \frac{1}{2} \begin{bmatrix} R^{(0,1)} \\ R^{(0,1)} \end{bmatrix}, \qquad \hat R^{(0,\ccT)} = \frac{1}{2} \begin{bmatrix} R^{(0,\ccT)} \\ R^{(0,\ccT)} \end{bmatrix} .
\end{equation}
For the agents in $\ccL_C$ the costs to be on an edge, sink or source is defined as before, i.e, for $\ell \in \ccL_C$ it is given by \eqref{eq:cost_trafic_example}.
Trucks are incentivized to use highways as much as possible by an increased cost for agents in $\ccL_T$ to be on small roads.
Thus, we define the modified cost matrix $\hat C_L \in \RR_{+}^{2L \times n}$ by
\begin{equation}
(\hat C_L)_{\ell i} = \begin{cases} ( C_L)_{\ell i}, & \text{ if }  \ell \in \ccL_C \\
 0.01,  & \text{ if } \ell \in \ccL_T, i \in \ccS^+ , \\
 0.1, &  \text{ if } \ell \in \ccL_T, i \in \ccH , \\
  0.7, &  \text{ if } \ell \in \ccL_T, i \in \ccE \setminus \ccH , \\
 0, & \text{ if } \ell \in \ccL_T, i \in \ccS^- . \end{cases}
\end{equation}
The rest of the problem is set up as before, and we solve it with Algorithm~\ref{alg:multi_commodity} and regularization parameter $\epsilon= 0.01$.
For each of the three sinks highlighted in Figure~\ref{fig:sthlm_map}, we consider the two associated commodities, and show the number of agents on the small roads and highways over time in Figure~\ref{fig:status_3comm_trucks_roads}.
As expected, the trucks avoid the small roads and mainly use the highways. In order to not exceed the capacity constraints on the highways, the cars are thus forced to the small roads.

\section{Conclusion.}

We have developed a novel framework for dynamic network flow problems, which is based on formulating the problem as a structured multi-marginal optimal transport problem.
Regularizing the problem with an entropy term opens up for efficiently finding an approximate solution.
By taking advantage of the graph-structure in the optimal transport formulations, we derived a scheme that is computationally highly efficient, as well as easy to implement.
Its competitiveness with state-of-the-art methods for network flow problems is experimentally illustrated in performance studies and on a traffic routing problem with a huge number of commodities.

% Appendix here
% Options are (1) APPENDIX (with or without general title) or 
%             (2) APPENDICES (if it has more than one unrelated sections)
% Outcomment the appropriate case if necessary
%
% \begin{APPENDIX}{<Title of the Appendix>}
% \end{APPENDIX}
%
%   or 
%
% \begin{APPENDICES}
% \section{<Title of Section A>}
% \section{<Title of Section B>}
% etc
% \end{APPENDICES}

% Acknowledgments here
%\section*{Acknowledgments.}
% Enter the text of acknowledgments here

% References here (outcomment the appropriate case) 

% CASE 1: BiBTeX used to constantly update the references 
%   (while the paper is being written).
%\bibliographystyle{informs2014} % outcomment this and next line in Case 1
%\bibliography{<your bib file(s)>} % if more than one, comma separated

% CASE 2: BiBTeX used to generate mypaper.bbl (to be further fine tuned)
%\input{mypaper.bbl} % outcomment this line in Case 2

%\bibliographystyle{plain}
\bibliographystyle{informs2014}
%\bibliography{BibJka,ref}
\bibliography{ref}

\end{document}